\documentclass[10pt]{amsart}
\usepackage{amsmath}
\usepackage{amsthm}
\usepackage{amsfonts}
\usepackage{amssymb}
\usepackage[all]{xy}

\newcommand{\Ps}{\mathbb{P}}
\newcommand{\Spg}{\overline{S_g}}

\newcommand{\ra}{\rightarrow}
\newcommand{\lra}{\longrightarrow}
\newcommand{\PNg}{\mathbb P_{N_g}}

\newcommand{\col}{\colon}

\renewcommand{\phi}{\varphi}
    
    \newtheorem{Lem}{Lemma}[section]
    \newtheorem{Prop}[Lem]{Proposition}
    \newtheorem*{thm1}{Theorem 1}
    \newtheorem*{thm2}{Theorem 2}
    \newtheorem*{thm3}{Theorem 3}

    \newtheorem{Thm}[Lem]{Theorem}

\theoremstyle{definition}

    \newtheorem{Def}[Lem]{Definition}
    \newtheorem{Exa}[Lem]{Example}
    \newtheorem{Rem}[Lem]{Remark}

    \DeclareMathOperator{\spa}{span}    
    \DeclareMathOperator{\first}{i}     
        
    \DeclareMathOperator{\Pic}{Pic}
    
    \DeclareMathOperator{\Hilb}{Hilb}

\begin{document}

\title{Spin curves over non-stable curves}

\author{Marco Pacini}

\address{Dipartimento di Matematica Guido Castelnuovo, Universit\`a Roma La Sapienza, P.le Aldo Moro 2, 00185 Roma, Italia}

\address{Instituto de Matematica Pura e Aplicada, Estrada D. Castorina 110, 22460-320 Rio de Janeiro, Brazil}

\address{Universidade Federal Fluminense, Rua M\'ario Santos Braga, s/n, Valonguinho, 24020-005 Niter\'oi RJ, Brazil }

\email{pacini@impa.br}

\thanks{The author was parcially supported by CNPq, Proc.151610/2005-3, and by Faperj, Proc. E-26/152-629/2005}

\begin{abstract}
Here we consider degenerations of stable spin curves for a fixed smoothing of a non-stable curve: we are able to give enumerative results and a description of limits of stable spin curves. We give a geometrically meaningful definition of spin curves over non-stable curves.
\end{abstract}

\maketitle

\section{Introduction}
 
The problem of constructing a compactification for the Picard scheme (or generalized Jacobian) of a singular algebraic curve has been studied by several authors. More generally, the same problem can be considered for families of curves. 

Several constructions have been carried out since Igusa's work \cite{I}, which gave a construction for nodal and irreducible curves. Constructions are known for families of geometrically integral curves, by Altman and Kleiman \cite{AK}, and geometrically connected, possibly reducible, nodal curves, by Oda and Seshadri \cite{OS}. A common approach to the problem is the use of the Geometric Invariant Theory. We recall in particular Caporaso's \cite{Capth} and Pandharipande's \cite{P} modular compactifications of the universal Picard variety over the moduli space of Deligne-Mumford stable curves. A different method was employed by Esteves \cite{E} to produce a compactification for a family of geometrically reduced and connected curves. 

\smallskip

On the other hand, one may be interested in distinguished subschemes of the Picard scheme. In the paper \cite{Corn}, Cornalba constructed a geometrically meaningful compactification $\Spg$ of the moduli space of theta characteristics of smooth curves of genus $g.$ The moduli space $\Spg$ is well-known as \emph{moduli space of stable spin curves} and is endowed with a natural finite morphism $\phi\col\Spg\lra\overline{M_g}$ onto the moduli space of Deligne-Mumford stable curves. As one can expect, the degree of $\phi$ is $2^{2g}$ and $\Spg$ is a disjoint union of two irreducible components, $\overline{S^+_g}$ and $\overline{S^-_g},$ whose restrictions over $M_g$ parametrize respectively even and odd theta characteristics on smooth curves. In particular, the degree of the restriction of $\phi$ to $\overline{S^-_g}$ is $N_g:=2^{g-1}(2^g-1).$ 

The fibers of $\phi$ over singular curves parametrize \emph{stable spin curves}, which are generalized theta-characteristics. The paper \cite{CC} provides an explicit combinatorial description of the boundary, parametrizing certain line bundles on quasistable curves having degree $1$ on exceptional components, i.e.  rational components intersecting the rest of the curve in exactly $2$ points.

More recently, in \cite{CapCasCorn}, the authors generalize the construction compactifying in the same spirit the moduli space of pairs $(C,L),$ $C$ a smooth curve and $L$ a $r$-th root of a fixed $N\in\Pic C.$

\smallskip

In this paper we will often make the following assumptions:

\begin{itemize}
\item[(1)]
a one-parameter family of projective curves has local complete intersection (l.c.i.) fibers which are Gorenstein, reduced, connected and canonical; 
\item[(2)]
a singular curve is irreducible with at most nodal, cuspidal and tacnodal singularities.
\end{itemize}

We construct a compactification of the moduli space of odd theta characteristics on the smooth fibers of a family of curves satisfying (1) and (2).  These assumptions allow us to find a rather explicit geometric description of degenerations of odd theta characteristics. Our method gives the possibility to reduce ourselves to results on Deligne-Mumford stable curves. Loosely speaking, this approach can be viewed as a ``Stable Reduction for polarized curves''.

\smallskip

Let us give more details. We say that a one-parameter family $f\col\mathcal W\ra B,$ with $B$ an affine and connected smooth curve, is a smoothing of a curve $W,$ if its general fiber is smooth and the fiber over a special point $0\in B$ is $W.$ Let $f\col\mathcal W\ra B$ be a smoothing of a singular curve $W$, where $\mathcal W\subset B\times\Ps^{g-1}$. Assume that $f$ satisfies $(1).$ Set $B^*:=B-0$ and consider the restricted family $\mathcal W^*\ra B^*.$ It is well-known that there exists a curve $S^-_{\omega_f^*},$ finite over $B^*,$ whose points parametrize odd theta characteristics of the fibers of $\mathcal W^*\ra B^*.$ Some natural questions arise:

\begin{itemize}
\item[(a)]
how can we get a compactification of $S^-_{\omega_f^*}$ over $B,$ reflecting the geometry of $W$?
\item[(b)]
are the corresponding boundary points independent of the chosen family $f\col\mathcal W\ra B$?
\item[(c)]
if the answers to (a), (b) are positive, can we give a geometric description of the boundary?
\end{itemize}

It is well-known that a smooth curve $C$ of genus $g$ has exactly $N_g=2^{g-1}(2^g-1)$ odd theta characteristics. If $C$ is general, any such line bundle $L$ satisfies $h^0(C,L)=1.$ Thus the canonical model of $C$ admits exactly one hyperplane $H_L$ cutting the double of the effective divisor of $|L|$. In this case, we say that $C$ is \emph{theta-generic} and that $H_L$ is a \emph{theta hyperplane} of $C.$ It follows that, if $C$ is a theta-generic curve, it comes with a configuration of theta hyperplanes $\theta(C),$ a point of $\Hilb_{N_g}(\Ps^{g-1})^\vee.$ 
Let $\Hilb_{g-1}^{p(x)}$ be the Hilbert scheme of curves of $\Ps^{g-1}$ having $p(x)=(2g-2)x-g+1$ as Hilbert polynomial and let $H_g$ be the irreducible component of $\Hilb_{g-1}^{p(x)}$ whose general point parametrizes a smooth canonical curve.  Consider the rational map:
\[
\SelectTips{cm}{11}
\begin{xy} <16pt,0pt>:
\xymatrix{
\theta\col H_g \UseTips\ar@{.>}[r] & \Hilb_{N_g}(\Ps^{g-1})^\vee
}
\end{xy}
\]
sending (a point parametrizing) a smooth theta-generic curve to its configuration of theta hyperplanes. If the smooth fibers of $f\col\mathcal W\ra B$ are theta-generic, the family of theta hyperplanes of $\mathcal W^*\ra B^*$ is isomorphic to $S^-_{\omega_f^*},$ hence its closure in $B\times (\Ps^{g-1})^\vee$ provides a compactification of $S^-_{\omega_f^*}.$ In this way, we can also consider limit theta hyperplanes on singular curves; a singular curve is \emph{theta-generic} if it admits a finite number of theta hyperplanes.

The following Theorem 1 answers question (b) for certain types of curves. For a proof of Theorem 1, see the proof of Proposition \ref{LCI} and Lemma \ref{zero-type}.

\begin{thm1} 
Let $W$ be a l.c.i. canonical curve. Then the following statements hold.
\begin{itemize}
\item[(i)]
Fix non negative integers $\tau, \gamma, \delta.$ If $W$ is general, irreducible  with $\tau$ tacnodes, $\gamma$ cusps and $\delta$ nodes, then it is theta-generic.
\item[(ii)]
If $W$ is theta-generic, then there exists a natural configuration of theta hyperplanes $\theta(W)\in \Hilb_{N_g}(\Ps^{g-1})^\vee$ such that, when $W$ is smooth, $\theta(W)$ is the ordinary configuration of theta hyperplanes.
\end{itemize}
\end{thm1} 

We are able to give an explicit description of $\theta (W)$ as follows. If $W$ is a general irreducible l.c.i. canonical curve with tacnodes, cusps and nodes, we say that a theta hyperplane of $W$ is \emph{of type $(i,j,k,h)$} if it contains $i$ tacnodes and $j$ tacnodal tangents of these $i$ tacnodes, $k$ cusps and $h$ nodes of $W.$ Denote by $t^j_{ikh}(W),$ for $j\le i,$ the number of theta hyperplanes of type $(i,j,k,h)$. Set $N^+_g:=2^{g-1}(2^g+1)$. 

The following Theorem 2 extends known results from \cite{Cap}. For a proof of Theorem 2, see the proof of Theorem \ref{hyper}.

\begin{thm2}
Let $g$ be an integer $g\ge 3.$ Fix non negative integers $\tau, \gamma, \delta.$ Let $W$ be a general irreducible l.c.i. canonical curve of genus $g$ with $\tau$ tacnodes, $\gamma$ cusps and $\delta$ nodes. Let $\widetilde{g}=g-\delta-\gamma-2\tau$ be the genus of the normalization of $W.$ 

If $j<i$ or $h\ne\delta,$ then: $$t^j_{ikh}(W)=2^{2\widetilde{g}+\tau-j+\delta-h-1}\binom{\tau}{i}\binom{i}{j}\binom{\delta}{h}\binom{\gamma}{k}.$$

If $i=j$ and $h=\delta,$ then: $$t^i_{ik\delta}(W)=\left\{\begin{array}{ll}\displaystyle   2^{\tau-i}\binom{\tau}{i}\binom{\gamma}{k} N_{\widetilde{g}}& \text{if \ }\tau -i+ \gamma -k\equiv 0 \;(2) \\ \\ 
\displaystyle 2^{\tau-i}\binom{\tau}{i}\binom{\gamma}{k}N^+_{\widetilde{g}}&\text{if \ }\tau-i+\gamma-k\equiv 1 \; (2);.\end{array}\right.$$\end{thm2}

If $W$ is singular, then $\theta (W)$ contains multiple hyperplanes. In  Theorem 3, we find the multiplicity of a limit theta hyperplane, as a multiplicative function of the singularities of $W.$ For a proof of Theorem 3, see the proof of Theorem \ref{mult}.

\begin{thm3}
Let $W$ be a general irreducible l.c.i. canonical curve with tacnodes and cusps. The multiplicity of a theta hyperplane of type $(i, j, k)$ is $4^{i-j}\;6^j\;3^k.$ 
\end{thm3}

The techniques used to prove Theorem 3 also lead to answer question (c). Let us  start with an example. Consider a general smoothing $\mathcal W\ra B$ of a projective irreducible canonical curve $W$ with one cusp. Modulo a base change we can assume that it admits a stable reduction over $B$, which we denote by $f\col\mathcal C\ra B.$ The central fiber $C$ of $\mathcal C$ is reducible. There exists a morphism from $\mathcal C$ to $\mathcal W$ given by $\mathcal N=\omega_f(D),$ a twist of the relative dualizing sheaf $\omega_f$ by a non-trivial divisor $D$ of $\mathcal C,$ supported on an irreducible component of $C.$ This morphism encodes the stable reduction of the polarized curve $(W, \mathcal O_W(1))$, suggesting a geometrically meaningful connection between limit theta characteristics on $W$ and square roots of the restriction $\mathcal N|_C$. We will explicitly describe this connection. We define a \emph{twisted spin curves} as a square roots of a twist of the dualizing sheaf of nodal curves. For example the square roots of $\mathcal N|_C$ are twisted spin curves. We will see that, if $W$ is as in Theorem 3, then the hyperplanes of $W$ correspond to suitable twisted spin curves of the curve, which is the stable reduction of any general smoothing of $W.$

\smallskip

In short, in Section \ref{sec2}, we give a review of moduli spaces of line bundles of curves. In Section \ref{sec3}, we introduce our compactification of $\overline{S}^-_{\omega_f^*}$ and we prove the existence of a well-defined configuration of theta hyperplanes for certain singular curves. In Section \ref{sec4} and Section \ref{sec5}, we give enumerative results of configurations of theta hyperplanes, describing their zero-dimensional scheme. In Section \ref{sec5.4}, we conclude with a definition of spin curves over non-stable curves.

\smallskip

We will use the following notation and terminology. We work over the field of complex numbers. A \emph{curve} is a connected projective curve which is Gorenstein and reduced. Let $\omega_W$ be the dualizing sheaf of a curve $W$. The genus of a curve is $g_W=h^0(W,\omega_W).$  We denote by $\text{Sing}(W)$ the set of the singularities of $W$. If $Z\subset W$ is a subcurve, set $Z^c:=\overline{W-Z}.$ A \emph{l.c.i.} curve is a local complete intersection curve. A curve with a \emph{cusp} or a \emph{tacnode} is a curve on a smooth surface with a singularity of curves of type $A_2$ or $A_3$, i.e. a planar singularity of a curve which locally analitically has equation $y^2=x^3$ or $y^2=x^4$.

A \emph{family of curves} is a proper and flat morphism $f\col\mathcal W\ra B$ whose fibers are curves. The fiber of a family $f\col\mathcal W\ra B$ over the point $b\in B$ is denoted by $W_b.$ If $0$ is a distinguished point of a scheme $B,$ we denote by $B^*:=B-0.$ A \emph{smoothing} of a curve $W$ is a family $f\col\mathcal W\ra B,$ where $B$ is a smooth, connected, affine curve of finite type, with a distinguished point $0\in B,$ such that $W_0$ is isomorphic to $W$ and $W_b$ is smooth for $b\in B^*.$ A \emph{general smoothing} is a smoothing with smooth total space. If $f\col\mathcal W\ra B$ is a smoothing of $W,$  we denote by $\mathcal W^*$ the restriction of $\mathcal W$ over $B^*.$ Similarly, if $\mathcal N\in\Pic\mathcal W,$ we denote by $\mathcal N^*:=\mathcal N|_{\mathcal W^*}.$

A \emph{stable} (\emph{semistable}) curve $C$ is a nodal curve such that every smooth rational subcurve of $C$ meets the rest of the curve in at least $3$ points ($2$ points). The \emph{dual graph} $\Gamma_X$ of a nodal curve $X$ is the usual graph, with the irreducible components of $X$ as vertices and the nodes of $X$ as edges. A curve $X$ is obtained from $C$ by \emph{blowing-up} a subset $\Delta$ of the set of the nodes of $C,$ if there is a morphism $\pi\col X\ra C$ such that, for every $n_i\in\Delta,$ $\pi^{-1}(n_i)=E_i\simeq\Ps^1$ and $\pi\col X-\cup_i E_i\ra C-\Delta$ is an isomorphism. For every $n_i\in\Delta,$ we call $E_i$ an \emph{exceptional component}. A \emph{quasistable} curve is a semistable curve, obtained by blowing-up a stable curve. 

A non-degenerate curve $W\subset\Ps^{g-1}$ of genus $g$ is \emph{canonical} if $\mathcal O_W(1)\simeq\omega_W.$

We set $N_g:=2^{g-1}(2^g-1)$ and $N^+_g:=2^{g-1}(2^g+1),$ respectively the numbers of odd and even theta characteristics of a smooth curve of genus $g.$

\section{Review of moduli of roots of line bundles of curves}\label{sec2}

In \cite{CapCasCorn}, the authors focused on the problem of giving a compactification of moduli spaces of roots of line bundles on smooth curves. The compactification is described in terms of limits square roots.

Let $C$ be a nodal curve and let $N\in\text{Pic}(C)$ be of degree divisible by $2.$
A triple $(X,L,\alpha),$ where $\pi\col X\ra C$ is a blow-up of $C,$ $L$ is a line bundle on $X$ and $\alpha$ is a homomorphism $\alpha\col L^{\otimes 2}\ra \pi^*(N),$ is a \emph{limit square root} of $(C,N)$ if:
\begin{itemize}
\item[(i)]
the restriction of $L$ to every exceptional component has degree $1;$
\item[(ii)]
the map $\alpha$ is an isomorphism at the points of $X$ not belonging to an exceptional component;
\item[(iii)]
for every exceptional component $E$ such that $E\cap E^c=\{p,q\}$ the orders of vanishing of $\alpha$ at $p$ and $q$ add up to $2.$ 
\end{itemize}

The curve $X$ is called the \emph{support} of the limit square root. If $C$ is stable, then a limit square root of $(C, \omega_C)$ is said to be a \emph{stable spin curve}.

If $X$ is a quasistable curve, we set $\widetilde{X}:=\overline{X-\cup E},$ where $E$ runs over the set of the exceptional components. We denote by $\Sigma_X$ the graph having the connected components of $\widetilde{X}$ as vertices and the exceptional components as edges. There exists a notion of isomorphism of limit square roots. By \cite[Lemma 2.1]{Corn}, two limit square roots $\xi=(X,L,\alpha)$ and 
$\xi'=(X,L',\alpha')$ are isomorphic if and only if the restrictions of $L$ and $L'$ to $\widetilde{X}$ are isomorphic. When there is no possibility of confusion, we denote by $\xi=(X,L,\alpha)$ both a limit square root and its isomorphism class.

Let $f\col\mathcal C\ra B$ be a family of nodal curves over a quasi-projective scheme $B$ and let $\mathcal N\in\Pic(\mathcal C)$ be of even relative degree. There exists a quasi-projective scheme $\overline{S}_f(\mathcal N),$ finite over $B$, which is a coarse moduli space, with respect to a suitable functor, of isomorphism classes of limit square roots of the restriction of $\mathcal N$ to the fibers of $f$. For more details,  we refer to \cite[Theorem 2.4.1.]{CapCasCorn}.
Let $C$ be a nodal curve and $N\in\Pic (C)$ of even degree. Denote by $\overline{S}_C(N)$ the zero-dimensional scheme $\overline{S}_{f_C}(N),$ where $f_C\col C\ra\{pt\}$ is the trivial family. In particular, $\overline{S}_C(N)$ is in bijection with the isomorphism classes of limit square roots of $(C, N)$. If $f:\mathcal C\ra B$ is a family of curves and $\mathcal N\in\Pic\mathcal C,$ then the fiber of $\overline{S}_f(\mathcal N)\ra B$ over $b\in B$ is $\overline{S}_{C_b}(\mathcal N|_{C_b}),$ as explained in \cite[Remark 2.4.3]{CapCasCorn}.

\begin{Def}\label{admgraph}
Fix a blow-up $\pi\col X\ra C$ of a stable curve $C.$ The graph $A_X$ \emph{associated to $X$} is the subgraph of the dual graph $\Gamma_C$ of $C,$ whose edges correspond to the set of nodes of $C,$ which are blown-up by $\pi.$
A subgraph $A$ of $\Gamma_C$ is \emph{admissible} if for every irreducible component $C_j$ of $C,$ whose corresponding vertex of $\Gamma_C$ is $v_j,$ then the number of edges of $A$ containing $v_j$ is congruent to $\text{deg}_{C_j}(N)$ modulo $2.$ 
\end{Def}

Recall that, by \cite[2.2]{CapCasCorn}, a subgraph $A$ of $\Gamma_C$ is the graph associated to a blow-up $X$ of $C$ such that $X$ is the support of some limit square root of $(C,N),$ if and only if $A$ is admissible. There are $2^{b_1(\Gamma_C)}$ admissible subgraphs of $\Gamma_C.$ 

Let $A_X$ be an admissible subgraph of $\Gamma_C$. Denote by $E_1,\dots,E_m$ the exceptional components and by $E_i\cap E^c_i=\{p_i,q_i\}.$ Consider the restriction $\widetilde{\pi}\col\widetilde{X}\ra C$ of the blow-up morphism. By the given definitions, the dual graph of $\widetilde{X}$ is $\overline{\Gamma_C-A_X}.$ If $g^\nu$ is the genus of the normalization $C^\nu$ of $C,$ then there are $2^{2g^\nu+b_1(\overline{\Gamma_C-A_X})}$ line bundles $\widetilde{L}\in\text{Pic}(\widetilde{X})$ such that:
\begin{equation*}
\textstyle
\widetilde{L}^{\otimes 2}=\widetilde{\pi}^*(N)(-\sum_{1\le i\le m}(p_i+q_i)).
\end{equation*}

Indeed we have $2^{2g^\nu}$ choices for the pull-back of $\widetilde{L}$ to $C^\nu$ and $2^{b_1(\overline{\Gamma_C-A_X})}$ gluings at nodes of $\widetilde{X}.$
A limit square root of $(C,N)$ supported on $X$ is given by gluing any $\widetilde{L}$ to $\mathcal O_{E_i}(1)$ for $i=1,\dots,m.$

By \cite[4.1]{CapCasCorn}, the multiplicity of $\overline{S}_C(N)$ in $\xi=(X,G,\alpha)$ is $2^{b_1(\Sigma_X)}.$

\begin{Exa}\label{banana}
Let $C=\cup_{0\le j\le N}C_j$ be a stable curve, with dual graph shown below.
\[
\begin{xy} <16pt,0pt>:
(0,0)*{\scriptstyle\bullet}="a"; 
(0,-1.5)*{\scriptstyle\bullet}="b";
(1.5,0)*{\scriptstyle\bullet}="d";
(-1.5,0)*{\scriptstyle\bullet}="f";
"a"+0;"b"+0**\crv{(-0.8,-0.8)};
"a"+0;"b"+0**\crv{(0.8,-0.8)}; 
"a"+0;"d"+0**\crv{(0.8,-0.8)};
"a"+0;"d"+0**\crv{(0.8,0.8)}; 
"a"+0;"f"+0**\crv{(-0.8,-0.8)};
"a"+0;"f"+0**\crv{(-0.8,0.8)}; 
"a"+(0,-3)*{\text{The dual graph of }C};
"b"+(0,-0.5)*{C_2};
"d"+(0.5,0)*{C_1}; 
"f"+(-0.7,0)*{C_j}; 
"a"+(0,0,6)*{C_0};
\end{xy}
\]

We describe $\overline{S}_C(\omega_C),$ the zero-dimensional scheme of stable spin curves of $C.$ Let $X\ra C$ be a blow-up of $C.$ Then $A_X$ is admissible if and only if for $1\le j\le N$ either the two edges connecting $C_0$ to $C_j$ appear in $A_X$ or none appear. Let $A_X$ be the admissible graph of the blow-up $X$ of $C$ at the first $r$ pairs of nodes. A stable spin curve is $(X, G, \alpha),$ where $G$ is a gluing of a square root of $X_0\cup X_{r+1}\cup\cdots\cup X_N,$ of a square root of $C_j$ for $1\le j\le r$ to $\mathcal O_E(1)$ for any exceptional component $E.$ Since $b_1(\Sigma_X)=r,$ we have that $(X, G, \alpha)$ has multiplicity $2^r$ in $\overline{S}_C(\omega_C).$

Fix a smoothing $\mathcal C$ of $C.$ Let $D$ be a Cartier divisor of $\mathcal C$ supported on $C_1,\dots,C_N.$ Pick $T:=\mathcal O_{\mathcal C}(D)\otimes \mathcal O_C.$ A limit square root of $(C, \omega_C(T))$ supported on $C$ is simply a square root of $\omega_C(T).$ Being $b_1(\Sigma_C)=0,$ the multiplicity of such limit square root in $\overline{S}_C(\omega_C(T))$ is $1.$
\end{Exa}

Let $X$ be a nodal curve and fix a smoothing $f\col\mathcal X\ra B$ of $X.$ A line bundle $T\in\text{Pic}(X)$ is said to be a \emph{twister} of $X,$ if $T\simeq\mathcal O_{\mathcal X}(D)|_X,$ where $D$ is a Cartier divisor of $\mathcal X$ supported on irreducible components of $X.$ When there is no possibility of confusion, we denote a twister of $X$ by $\mathcal O_f(D).$

\begin{Def}
Let $L\in\text{Pic}(X)$ be a line bundle such that $L^{\otimes 2}\simeq\omega_X\otimes T$, where $T=\mathcal O_{\mathcal X}(D)\otimes\mathcal O_X$ is a twister of $X.$ Then $(X,L)$ is a \emph{D-twisted spin curve} or simply a \emph{twisted spin curve}. 
\end{Def}

A stable spin curve supported on a stable curve is a $0$-twisted spin curve. Notice that a twisted spin curve $(X, L)$ can be seen as a limit square root of $(X,\omega_X\otimes T).$

\subsection{The sections of a stable spin curve}\label{sec2.1} Let $C$ be a stable curve and let $\xi=(X,G,\alpha)$ be a stable spin curve of $C,$ supported on a blow-up $\pi\col X\ra C$ of $C.$ Let $\mathcal{E}(X)$ be the set of the exceptional components. Pick $\widetilde{X}=\overline{X-\cup_{E\in\mathcal{E}(X)}E}.$
The line bundle $G$ is obtained by gluing theta characteristics of the connected components of $\widetilde{X}$ to $\mathcal O_E(1)$ for every $E\in\mathcal{E}(X).$ Let $Z$ be a connected component of $\widetilde{X}.$ Since $G|_E=\mathcal O_E(1)$ for $E\in\mathcal{E}(X),$ a non-trivial section of $G|_Z$ uniquely extends to a section of $G$ vanishing on the other connected components of $\widetilde{X}.$ Thus:
\begin{equation}\label{decompos}
H^0(X,G)=\oplus H^0(Z,G|_Z)\;\;Z\text{ connected}.
\end{equation}

If $G$ is odd on $d_G$ connected components of $\widetilde{X},$ then $G$ is odd if and only if $d_G\equiv 1\;(2).$

\subsection{Smoothing line bundles and sections}\label{sec2.2} 
Let $f\col\mathcal W\ra B$ be a smoothing of a singular curve $W$ with nodes, cusps and tacnodes and let $\mathcal N\in\Pic(\mathcal W).$ 

Let $L\in\Pic(W)$ be endowed with an isomorphism $\iota_0\col L^{\otimes 2}\ra
\mathcal N\otimes\mathcal O_W.$ Then, up to shrinking $B$ to a complex neighbourhood of $0$, there exists a line bundle $\mathcal L\in\Pic\mathcal W$ extending $L$ and an isomorphism $\iota\col\mathcal L^{\otimes 2}\ra
\mathcal N$ extending $\iota_0.$ Moreover, if $(\mathcal L',i')$ is another extension of $(L,\iota_0),$ then there is an isomorphism $\chi\col\mathcal L\ra\mathcal L',$ restricting to the identity, and with $\iota=\iota'\circ\chi^{\otimes 2}.$ 

The previous statement is a simple variation of \cite[Remark 3.0.6.]{CapCasCorn}, whose proof appeared in an early version of the paper.

Assume that $h^0(W_b,\mathcal N|_{W_b})=d$ for $b\in B^*.$ Consider $f_*\mathcal N$ as a vector bundle. The \emph{$\mathcal N$-smoothable sections of $\mathcal N|_W$} are the sections of the $d$-dimensional subspace of $H^0(W,\mathcal N|_W)$ which is the fiber of $f_*\mathcal N$ over $0$.

\section{The projective setup of theta hyperplanes}\label{sec3}

A genus $g$ canonical smooth curve $W\subset\Ps^{g-1}$ has $N_g:=2^{g-1}(2^g-1)$ odd theta characteristics. If $W$ is general, then any odd theta characteristic $L$ has $h^0(W,L)=1$. Thus a general smooth canonical curve has exactly $N_g$ hyperplanes, called \emph{theta hyperplanes}, cutting the double of a semicanonical divisor. We call $W$ \emph{theta-generic}. We collect these hyperplanes, in a configuration $\theta(W)\in\PNg:=\Hilb_{N_g}(\Ps^{g-1})^\vee.$ The theta hyperplanes were introduced in \cite{CapSerbi}, \cite{Cap} and \cite{CapSerth}.
 
We define configurations of theta hyperplanes for singular curves as follows. Let $\Hilb_{g-1}^{p(x)}$ be the Hilbert scheme of curves of $\Ps^{g-1}$ having $p(x)=(2g-2)x-g+1$ as Hilbert polynomial and let $H_g$ be the irreducible component of $\Hilb_{g-1}^{p(x)}$ whose general point parametrizes a smooth canonical curve. 
If $h\in H_g$, let $W_h$ be the curve represented by $h.$   
Consider the map:
\[
\SelectTips{cm}{11}
\begin{xy} <16pt,0pt>:
\xymatrix{
\theta\col H_g \UseTips\ar@{.>}[r] & \PNg
}
\end{xy}
\]

such that $\theta(h)=\theta(W_h),$ for every smooth theta-generic canonical curve $W_h.$ Let $f:\mathcal W\ra B$ be a smoothing of a canonical curve $W\subset\Ps^{g-1}$ to theta-generic canonical curves, where $\mathcal W\subset B\times\Ps^{g-1}$. Let $\gamma_f\col B^*\ra H_g$ be the morphism associated to the restricted family $\mathcal W^*\rightarrow B^*.$ The map $\theta$ is defined on the image of $\gamma_f$. Now, $B$ is smooth and $\PNg$ is projective, then the map:
\[
\SelectTips{cm}{11}
\begin{xy} <16pt,0pt>:
\xymatrix{
\theta\circ\gamma_f\col B\UseTips\ar[r] & \PNg
}
\end{xy}
\]

is defined over $0\in B$. We set $\theta_f(W):=\theta\circ\gamma_f(0)$. We can see $\theta_f(W)$ also as a not necessarily reduced hypersurface of degree $N_g$ in $\Ps^{g-1}$, all of whose irreducible components are hyperplanes. Notice that a priori the configuration $\theta_f(W)$ depends on $f.$ Furthermore, we can consider the $B$-curve:
\begin{equation}\label{newcompdef}
\SelectTips{cm}{11}
\begin{xy} <16pt,0pt>:
\xymatrix{
J_{\mathcal W}\UseTips\ar[r] & B
}
\end{xy}
\end{equation}

which is the closure in $B\times(\Ps^{g-1})^\vee$ of the incidence correspondence:
$$\{(b ,H)\;|\;H\subset\theta(W_b)\;,\;b\ne 0\}.$$

\begin{Def}
We say that the hyperplanes of the fiber of the morphism (\ref{newcompdef}) over $0\in B$ are \emph{theta hyperplanes of $W$.} We say that $W$ is \emph{theta-generic} if it has a finite number of theta hyperplanes, arising from smoothings to theta-generic curves. 
\end{Def}

\begin{Prop}\label{LCI}
Let $W$ be a theta-generic l.c.i. canonical curve. If $f$ and $f'$ are two smoothings of $W$ to theta-generic canonical curves, then $\theta_f(W)=\theta_{f'}(W).$ 
 In particular, there is a natural configuration of theta hyperplanes $\theta(W)$ such that, when $W$ is smooth, $\theta(W)$ is the ordinary configuration of theta hyperplanes.
\end{Prop}

\begin{proof}
Assume that $H_g$ is smooth at the point parametrizing $W$. Consider  $H^{sm}_g$ and let $U\subset H^{sm}_g$ be the open set corresponding to theta-generic smooth curves. Let $\Gamma$ be the closure in $H^{sm}_g\times \PNg$ of the incidence variety:
$$\Gamma_U:=\{(h,\theta(W_h))|h\in U\}\subset H^{sm}_g\times \PNg.$$

Let $\rho$ be the projection $\rho\col\Gamma\ra H^{sm}_g.$ Since $\rho$ is bijective on $\Gamma_U,$ it is a birational morphism. Since $U$ is irreducible, also $\Gamma$ is irreducible. Thus, by the Zariski Main Theorem, the fibers of 
$\rho$ are connected. The curve $W$ is theta-generic, hence the fiber of $\rho$ over the point parametrizing $W$ is finite, hence it consists of one element and we are done. We show that, if $h\in H_g$ parametrizes $W,$ then $H_g$ is smooth at $h.$ If $W$ is l.c.i., it is enough to show that $h^1(N_{W/\mathbb{P}^{g-1}})=0$. Consider the sequence:
\[
\SelectTips{cm}{11}
\begin{xy} <16pt,0pt>:
\xymatrix{
0 \UseTips\ar[r] & \mathcal I_W/\mathcal I^2_W \ar[r] & \Omega^1_{\mathbb P^{g-1}|_W}\ar[r]& \Omega^1_W \ar[r] & 0
}
\end{xy}
\]

which is exact, since $W$ is l.c.i.. By taking $\mathcal Hom_{\mathcal O_W}(-,\mathcal O_W),$ we get:
\[
\SelectTips{cm}{11}
\begin{xy} <16pt,0pt>:
\xymatrix{
0 \UseTips\ar[r] & \mathcal Hom_{\mathcal O_W}(\Omega^1_W,\mathcal O_W)  \ar[r] & \mathcal T_{\mathbb{P}^{g-1}}|_W \ar[r]& N_{W/\mathbb{P}^{g-1}}
\ar[r]^{\,\,\,\,\,\,\,\,\,\alpha} &}
\end{xy}
\]
\[
\SelectTips{cm}{11}
\begin{xy} <16pt,0pt>:
\xymatrix{\ar[r]^{\alpha\hskip1.3cm} & \mathcal Ext^1_{\mathcal O_W}(\Omega^1_{W},\mathcal O_{W}) \ar[r] & 0.
}
\end{xy}
\]
Let $\mathcal N'_W$ be the kernel of $\alpha.$ By the sequences in cohomology, we get:
\[
\SelectTips{cm}{11}
\begin{xy} <16pt,0pt>:
\xymatrix{
H^1(W,\mathcal T_{\Ps^{g-1}|_W}) \UseTips\ar[r] & H^1(W,\mathcal N'_W) 
\ar[r] & 0
}
\end{xy}
\]
\[
\SelectTips{cm}{11}
\begin{xy} <16pt,0pt>:
\xymatrix{
H^1(W,\mathcal N'_W)\UseTips\ar[r] & H^1(W,N_{W/\Ps^{g-1}})\ar[r] & 0.
}
\end{xy}
\] 

Thus, if $h^1(W,\mathcal T_{\mathbb{P}^{g-1}|_W})=0,$ then $h^1(W,\mathcal N'_W)=0$ and also $h^1(W,N_{W/\Ps^{g-1}})=0.$ By the Euler sequence of $\Ps^{g-1},$ restricted to $W,$ we have:
\[
\SelectTips{cm}{11}
\begin{xy} <16pt,0pt>:
\xymatrix{
H^1(W, \mathcal O_W) \UseTips\ar[r] & H^1(W, \mathcal O_W(1))\otimes H^0 (W, \mathcal O_W(1))^\vee\ar[r] &
}
\end{xy}
\] 
\[
\SelectTips{cm}{11}
\begin{xy} <16pt,0pt>:
\xymatrix{\ar[r] &H^1(W, \mathcal T_{\Ps^{g-1}}|_W)\ar[r] & 0.
}
\end{xy}
\] 

Since $\mathcal O_W(1)\simeq\omega_W,$ dualizing the last sequence, we get:
\[
\SelectTips{cm}{11}
\begin{xy} <16pt,0pt>:
\xymatrix{
0 \UseTips\ar[r] & H^1(W,\mathcal T_{\Ps^{g-1}}|_W)^\vee\ar[r] & H^0(W,\mathcal O_W)\otimes H^0(W,\omega_W) \ar[r]^{\,\,\,\,\,\,\,\,\,\,\,\,\,\,\,\,\,\,\,\,\beta} & H^0(W, \omega_W).
}
\end{xy}
\] 

Being $\beta$ injective, then $h^1(W,T_{\Ps^{g-1}}|_W)=0$ and we are done. 
\end{proof}

We call $\theta(W)$ \emph{the configuration of theta hyperplanes of $W$}. Its associated zero-dimensional scheme is \emph{the zero-dimensional scheme of theta hyperplanes of $W$}.

\section{Enumerative results}\label{sec4}

In this section, we write down formulas for the reduced zero-dimensional scheme of theta hyperplanes for singular curves. In \cite[Proposition 1, Proposition 4]{Cap} one can find formulas for nodes and cusps, which we generalize also for tacnodes.

\begin{Def}
Let $W$ be a genus $g$ irreducible curve with $\tau$ tacnodes, $\gamma$ cusps and $\delta$ nodes and let $\nu\col W^\nu\ra W$ be its normalization. We say that $W$ is \emph{general} if $(W^\nu, \nu^{-1}(\text{Sing(W)}))$ is general in $M_{\widetilde{g},n}$, where $n=2\delta+\gamma+2\tau$ and $\widetilde{g}=g-n+\delta$.  
\end{Def}

\begin{Def}
Let $W\subset\Ps^{g-1}$ be an irreducible general canonical curve with tacnodes, cusps and nodes. A hyperplane of $\Ps^{g-1}$ is \emph{of type $(i,j,k,h)$ with respect to $W$} if it contains $i$ tacnodes, $j$ tacnodal tangents of these $i$ tacnodes, $k$ cusps and $h$ nodes. A hyperplane of type $(0,0,0,0)$  simply \emph{a hyperplane of type 0}. We denote by $t^j_{ikh}(W)$ the number (if it is finite) of theta hyperplanes of $W$ of type $(i,j,k,h)$.  
\end{Def}

Notice that the divisor cut on $W$ by a theta hyperplane is the limit of a family of divisors on smooth curves such that each point of the support appears with an even coefficient. Consider the projection of a canonical integral curve $W\subset\Ps^{g-1}$ with tacnodes, cusps and nodes, from a singular point $s\in W$. Call $W'\subset\Ps^{g-2}$ the projected curve, which is a canonical curve. 
Let $H\subset\Ps^{g-1}$ be a hyperplane. If $s\in H$, call $H'\subset\Ps^{g-2}$ the projected hyperplane. Lemma \ref{cut} will implies that $H$ is a theta hyperplane of $W$ if and only if $H'$ is a theta hyperplane of $W'$. Notice that if $s\in W$ is a tacnode, then $s$ projects to a node $n\in W'$.

\begin{Lem}\label{trivial}
Let $W'$ be an irreducible general curve of genus $g(W')\ge 2$ with tacnodes, nodes and cusps and $p\in W'$ be a general point. Let $(W')^{sm}$ be the smooth locus of $W'$. If $g(W')=2$, set $W^0=(W')^{sm}-\{p,p'\}$, where 
$p+p'\in |\omega_{W'}|$. If $g(W')>2$, set $W^0=(W')^{sm}-p$. Then there is an \'etale morphism $b\col \widetilde{W}^0\ra W^0$ and a family of curves $\psi\col\mathcal W\ra \widetilde{W}^0$, with a section $s\col \widetilde{W}^0\ra\mathcal W$, such that:
\begin{itemize}
\item[(i)]
if $q'\in \widetilde{W}^0$, then $s(q')$ is a tacnode of $\psi^{-1}(q')$;
\item[(ii)]
the normalization of $\mathcal W$ at $s(\widetilde{W}^0)$ is $\pi\col W'\times \widetilde{W}^0\ra\mathcal W$;
\item[(iii)]
if $q'\in\widetilde{W}^0$, then $\pi^{-1}(s(q'))=\{p,b(q')\}$.
\end{itemize}
\end{Lem}

\begin{proof}
Set $g=g(W')+2$ and $N=\omega_{W'}(2p+2q)$, for $q\in W^0$. Since $h^0(N)=g+1$ and $\deg N=2g-2$, the linear system $|N|$ embeds $W'$ as a projective curve $W'_q\subset\Ps^g$. Abusing notation, we can see $p,q$ as points of $W'_q$. Consider a point $r$ of the line $\overline{pq}$, with $r\ne p,q$ and the projection $\pi_r$ from $r$. Call $X\subset\Ps^{g-1}$ the image of 
$W'_q$ via $\pi_r$ and $t=\pi_r(p)=\pi_r(q)$. Denote by $T_p$ the tangent of $W'_q$ at $p$ and by $T_q$ the tangent of $W'_q$ at $q$. Notice that $T_p$ does not contain $q$, because $h^0(N(-2p-q))=h^0(N(-2p))-1$, and similarly $T_q$ does not contain $p$.

We show that $\pi_r$ restricts to an isomorphism between $W'_q-\{p, q\}$ and $X-t$. In fact, consider the projection $\pi_t$ of $X$ from $t$. The image $Y\subset\Ps^{g-2}$ of $X$ via $\pi_t$ is the projection of $W'_q$ from $\overline{pq}$, hence it is the image of the morphism $W'\ra\Ps^{g-2}$ given by $|N(-p-q)|=|\omega_{W'}(p+q)|\simeq\Ps^{g-2}$. Pick the genus $g-1$ curve $\overline{Y}$ obtained from $W'$ by the nodal identification of $p$ and $q$. Let $n$ be the new node. Now, $\overline{Y}$ is not hyperelliptic, because $W'$ is general and $p+q\notin|\omega_{W'}|$ when $g(W')=2$, hence $\omega_{\overline{Y}}$ is very ample. The canonical model of $\overline{Y}$ is the image of the morphism $W'\ra\Ps^{g-2}$ given by $|\omega_{W'}(p+q)|$. Thus, 
$Y$ is the canonical image of $\overline{Y}$ and $t$ projects to $n$. In this way, $\pi_t\circ\pi_r$ restricts to an isomorphism between $W'_q-\{p, q\}$ and $Y-n$ and $\pi_r$ restricts to an isomorphism between $W'_q-\{p, q\}$ and $X-t$.

Observe that $\text{dim}H^0(N(-2p-2q))=g-2$, i.e. it is a codimension 3 subspace of $H^0(N)$ instead of a codimension 4 subspace, hence $T_p$ and $T_q$ meet in a point. Set $T_t=\pi_r(T_p)=\pi_r(T_q)\subset\Ps^{g-1}$. We show that $t$ is a double point. By the proof of \cite[Proposition (6.1)]{EGK}, it is enough to show the existence of a hyperplane of $\Ps^{g-1}$ intersecting $X$ at $t$ with multiplicity 2. Let $H\subset\Ps^{g-1}$ be a hyperplane containing $t$ and not cotaining $T_t$. Assume that $H$ intersects $X$ at $t$ with multiplicity $m$. If $m\ge 3$, then $H$ is the projection via $\pi_r$ of a hyperplane $\overline{H}\subset\Ps^g$ such that the multiplicities of the intersection of $\overline{H}$ and $W'_q$ at $p$ and of the intersection of $\overline{H}$ and $W'_q$ at $q$ sum-up to $m$. Hence $\overline{H}$ contains either $T_p$ or $T_q$ and $H$ contains $T_t$, a contradiction. Thus $m=2$ and $p$ is a double point. Notice that $t$ is not unibranch, because $W'_q$ is the normalization of $X$ at $t$ with $p,q$ lying over $t$. Furthermore, $T_t$ is the common tangent of the two branches of the sigularity. In this way, $t$ is an $A_{2k}$ singular point, for $k\ge 2$, i.e. if 
$(x,y)$ is an analytic coordinate system at $t$ of a smooth surface containing 
$X$, then the equation of $X$ is $y²-x^{2k}=0$, for $k\ge 2$.

Assume that $k>2$. Set $g_k=g+k-2$, $X_k=X$ and $t_k=t$. Now, $X_k$ has genus $g_k$ and it is not hyperelliptic, because $W'$ is general and $p+q\notin|\omega_{W'}|$ when $g(W')=2$. Let $\nu_k\col X_{k-1}\ra X_k$ be the partial normalization of $X_k$ at $t_k$, i.e. $X_{k-1}$ has an $A_{2k-2}$ singular point $t_{k-1}$ lying over $t_k$. Notice that $X_{k-1}$ is not hyperelliptic and has genus $g_{k-1}$. Consider the canonical model $W_k\subset\Ps^{g_k-1}$ of $X_k$. Then $W_k$ has degree $2g_k-2$. Consider the projection $\pi_k$ of $W_k$ from $t_k$ and call $W_{k-1}\subset\Ps^{g_k-2}$ the image of $\pi_k$. The map $\pi_k$ is given by the subspace of $H^0(\omega_{W_k})$ of dimension $g_{k-1}=g_k-1$ corresponding to the hyperplanes of $\Ps^{g_k-1}$ containing $t_k$. Let $H\subset\Ps^{g_k-1}$ be a general hyperplane containing $t_k$, cutting $W_k$ in $t_k$ with multiplicity 2. Let $D$ be the divisor given by the intersection of $H$ with the smooth locus of $W_k$. If $L=\mathcal O_{X_{k-1}}(\nu_k^*D)$, then  the composed map $\pi_k\circ\nu_k$ is given by a subspace of $H^0(L)$ of dimension $g_{k-1}$.  Now, $\deg L=2g_k-4=2g_{k-1}-2=\deg\omega_{X_{k-1}}$ and $h^0(L)\ge g_{k-1}$, thus, by Riemann-Roch, $h^0(\omega_{X_{k-1}}\otimes L^{-1})=h^0(L)-g_{k-1}+1\ge 1$. Since $X_{k-1}$ is irreducible, this implies $L\simeq\omega_{X_{k-1}}$. In particular, $\pi_k\circ\nu_k$ is given by $H^0(\omega_{X_{k-1}})$ and $W_{k-1}$ is the canonical model of $X_{k-1}$. Iterating the reasoning, we get projections $\pi_k,\pi_{k-1},\dots\pi_3$ such that the image $W_2$ of $\pi_3$ is a curve with a tacnode $t_2$ such that the normalization of $W_2$ at $t_2$ is $W'\ra W_2$, with $p,q$ lying over $t_2$.

Now, pick the trivial family $f\col W'\times W^0\ra W^0$. Consider the Cartier divisors $D=p\times W^0\subset W'\times W^0$ and $\Delta\subset W'\times W^0$, where $\Delta$ is the restriction of the diagonal of $W'\times W'$. 
Set $\mathcal N=\omega_f(2D+2\Delta)$. If $q\in W^0$, we have $\mathcal N|_{f^{-1}(q)}=\omega_{W'}(2p+2q)$. Embed the family $f$ as a family of projective curves of $\Ps^g$, via the relative linear system of $\mathcal N$. As before, call $W'_q$ the image of the fiber $f^{-1}(q)$. Consider: $$\mathcal S=\Ps(f_*\mathcal N/f_*\mathcal N(-D-\Delta))-(\Ps(f_*\mathcal N/f_*\mathcal N(-D))\cup\Ps(f_*\mathcal N/f_*\mathcal N(-\Delta)))$$ and the natural map $\theta\col\mathcal S\ra W^0$. The fiber of $\theta$ over $q\in W^0$ parametrizes the codimension 1 subspaces $V$ of $H^0(\mathcal N|_{f^{-1}(q)})$ such that: $$H^0(\mathcal N|_{f^{-1}(q)}(-p-q))\subset V\notin \{H^0(\mathcal N|_{f^{-1}(q)}(-p)), H^0(\mathcal N|_{f^{-1}(q)}(-q))\}.$$ Any such $V$ gives a projection of $W'_q$ from a point $r\in\overline{pq}$, $r\ne p,q$. 
Now, by \cite[IV-4, 17.16.3]{G}, there exists an \'etale morphism $b\col\widetilde{W}^0\ra W^0$ such that the morphism $\mathcal S\times_{W^0}\widetilde{W}^0\ra \widetilde{W}^0$ has a section.  Thus, for every $q'\in\widetilde{W}^0$ such that $b(q')=q$, we get the choice of a point $r\in\overline{pq}$, $r\ne p,q$, from which we can project $W'_q$. In this way, we get a family of curves $\psi_k\col\mathcal W_k\ra\widetilde{W}^0$, with a section $s_k\col\widetilde{W}^0\ra\mathcal W_k$ such that $s_k(q')$ is an $A_{2k}$ singular point of $\psi_k^{-1}(q')$, for some $k\ge 2$. Note that $k$ does not depend on $q'$, because the genus of $\psi_k^{-1}(q')$ does not. The normalization of $\mathcal W_k$ at $s_k(\widetilde{W}^0)$ is $\pi_k\col W'\times\widetilde{W}^0\ra\mathcal W_k$ with $\pi_k^{-1}(s_k(q'))=\{p,b(q')\}$.  If $k>2$, take the canonical model of the family, given by $\omega_{\psi_k}$, and the fiberwise projection from $s_k(q')$. We get a family $\psi_{k-1}\col\mathcal W_{k-1}\ra\widetilde{W}^0$, with a section $s_{k-1}\col\widetilde{W}^0\ra\mathcal W_{k-1}$, having the  properties of $\psi_k\col\mathcal W_k\ra\widetilde{W}^0$, with the exception that $s_{k-1}(q')$ is an $A_{2k-2}$ singular point. Iterating the reasoning, we get a family of curves $\psi_2\col\mathcal W_2\ra\widetilde{W}^0$ as required. 
\end{proof}

\begin{Lem}\label{zero-type}
Let $W$ be an irreducible general curve with tacnodes, cusps and nodes. Then the following statements hold.
\begin{itemize}
\item[(i)]
If $R$ is a theta characteristic of $W,$ then $h^0(R)\le 1$ and a section of $R$ does not vanish on singular points of $W$; 
\item[(ii)]
If $W$ is canonical, it is theta-generic and there exists a set bijection:
\[
\SelectTips{cm}{11}
\begin{xy} <16pt,0pt>:
\xymatrix{
\{\text{theta hyperplanes of type 0 of }W\}\UseTips\ar[r] & \{\text{odd theta characteristics of }W\}.
}
\end{xy}
\]
\end{itemize}
\end{Lem}

\begin{proof}
(i) Let $g$ be the genus of $W$. We argue by induction on the number of singularities of $W$. The statement is clear if $W$ is smooth, and also if $g\le 1$. Let $g\ge 2$ and let $R$ be a theta characteristic of $W$. 

If $W$ has a node $n,$ pick its normalization $\pi\col W'\ra W$ at $n$. Set $\pi^{-1}(n)=\{p, q\}.$ If a section of $R$ vanishes on $n$, then $h^0(\pi^*R(-p-q))\ne 0$. The genus of $W'$ is $g-1$ and $(\pi^*R)^{\otimes 2}=\pi^*(\omega_W)=\omega_{W'}(p+q)$. Since $\deg \pi^*R=g-1$, by Riemann-Roch $h^0(R)\le h^0(\pi^*R)=1+h^0(\pi^*R(-p-q))$. Thus, we are done if we show that 
$h^0(\pi^*R(-p-q))=0$. To show that $h^0(\pi^*R(-p-q))=0$, let $W^{sm}$ be the smooth locus of $W'$. Pick $\mathcal W=W'\times W^{sm}$ and the Cartier divisors $\Delta\subset\mathcal W$, where $\Delta$ is the restriction of the diagonal of $W'\times W'$, and $D=q\times W^{sm}\subset\mathcal W$. Let $f\col \mathcal W \ra W^{sm}$ be the second projection and $\mathcal M=\omega_f(\Delta+D)$. We have 
$\mathcal M|_{f^{-1}(p)}=\omega_{f^{-1}(p)}(p+q)$ and $\mathcal M|_{f^{-1}(q)}=\omega_{f^{-1}(q)}(2q)$. A square root of $\mathcal M|_{f^{-1}(q)}$ is 
$L(q)$, where $L$ is a theta characteristic of $f^{-1}(q)=W'$. Now, $q$ is general and by induction $h^0(L)\le 1$, then $h^0(L(-q))=0$. Let $\mathcal R$ be the line bunde smoothing $L(q)$ on a complex neighbourhood of 
$q\in W^{sm}$ such that $\mathcal R^{\otimes 2}=\mathcal M$, as in Section \ref{sec2.2}. If $p$ is a general point contained in such a complex neighbourhood and if a section of the square root $\mathcal R|_{f^{-1}(p)}$ of $\omega_{f^{-1}(p)}(p+q)$ vanishes in $p$ and $q$, then a section of $L(q)$ vanishes two times in $q$, i.e. $h^0(L(-q))\ne 0$, a contradiction. Now, the family $f$ is trivial, then the number of the square roots of the restriction of $\mathcal M$ to the fibers of $f$ is constant. Thus, arguing in the same way for all the square root of $\mathcal M|_{f^{-1}(q)}$, we get that the sections of any square root of $\omega_{f^{-1}(p)}(p+q)$ do not vanish on $p$ and $q$.

If $W$ has a cusp $c$, pick its normalization $\pi\col W'\ra W$ at $c.$ Set $\pi^{-1}(c)=p$. The genus of $W'$ is $g-1.$  Now, $(\pi^*R)^{\otimes 2}=\omega_{W'}(2p)$, then $\pi^*R=L(p),$ for a theta characteristic $L$ of $W'$. By induction, $h^0(L)\le 1$ and by Riemann-Roch, $h^0(\pi^* R)=1+h^0(\omega_{W'}\otimes(\pi^* R)^{-1})=1+h^0(L(-p)).$ Since $p$ is general, it follows that $h^0(L(-p))=0$ and then $h^0(R)\le h^0(\pi^*R)=1.$ Now, assume that 
$h^0(R)=1$ and that a section of $R$ vanishes on $c$. Thus $h^0(\pi^*R(-p))=h^0(L)\ne 0$. From the proof of \cite[Theorem (2.22)]{Ha} and \cite[Section 2e (a)]{Ha}, we have $h^0(R)\equiv h^0(L)+1\;(2)$, hence $1\ge h^0(L)\equiv 0\;(2)$. Thus $h^0(L)=0$, which is a contradiction.

If $W$ has a tacnodes $t$, pick its normalization $\pi\col W'\ra W$ at $t.$ Set $\pi^{-1}(t)=\{p, q\}.$ The genus of $W'$ is $g-2.$ Now, $(\pi^*R)^{\otimes 2}=\omega_{W'}(2p+2q)$, then $\pi^*R=L(p+q)$ for a theta characteristic $L$ of $W'$. By induction $h^0(L)\le 1.$ Since $p$ and $q$ are general, we have $h^0(L(-p-q))=0$, hence $h^0(R)\le h^0(\pi^* R)=2+h^0(L(-p-q))=2$. We are done if $h^0(R)=0$ for an even theta characteristic $R$. In fact, in this case, if $R$ is an odd theta characteristic, we have $h^0(R)=1$ and if the section of $R$ vanishes on $t$, then $h^0(\pi^*R(-p-q))=h^0(L)\ne 0$. But, from the proof of \cite[Theorem (2.22)]{Ha} and \cite[Section 2e (a)]{Ha}, we have $h^0(R)\equiv h^0(L)+1\;(2)$, hence $1\ge h^0(L)\equiv 0\;(2)$. Thus $h^0(L)=0$, which is a contradiction.

We show that $h^0(R)=0$ for an even theta characteristic $R$. This is well-known if $g\le 3$. Let $g>3$. Let $\psi\col\mathcal W\ra \widetilde{W}^0$ be the family of Lemma \ref{trivial}, associated to $W'$ and $p$, where $\widetilde{W}^0\ra W^0$ is \'etale and $W'\times\widetilde{W}^0\ra\mathcal W$ is the normalization at the distinguished tacnodes of the fibers. It is well-known that there exists an \'etale base change $Z\ra \widetilde{W}^0$ such that the relative degree $g-1$ Picard variety $J^{g-1}_{\mathcal Z/Z}$ of $\mathcal Z=Z\times_{\widetilde{W}^0}\mathcal W\ra Z$ has a universal object, i.e. a universal line bundle $\mathcal M$ over $J^{g-1}_{\mathcal Z/Z}\times_Z\mathcal Z$. Denote by $b_2\col Z\ra\widetilde{W}^0\ra W^0$ the composition of the \'etale morphisms. Let $\Theta_{\mathcal Z}\subset J^{g-1}_{\mathcal Z/Z}$ be the locus corresponding to theta characteristics and $b_1\col\Theta_{\mathcal Z}\ra Z$ the corresponding finite \'etale morphism. Set $\mathcal R=\mathcal M|_{\Theta_{\mathcal Z}\times_Z\mathcal Z}$ and consider the commutative diagram:
\[
\SelectTips{cm}{11}
\begin{xy} <16pt,0pt>:
\xymatrix{
\mathcal L \ar[r]& W' \times \Theta_{\mathcal Z}\ar[d]_{\rho_1}\UseTips\ar[r]^{\phi_1}& W'\times Z\ar[r]^{\phi_2}\ar[d]^{\pi}& W'\times W^0\ar[dd]\\
\mathcal R\ar[r] & \Theta_{\mathcal Z}\times_Z\mathcal Z \ar[d]_{\rho_2}\UseTips\ar[r]& \mathcal Z \ar[d]^{\eta} &  \\
 & \Theta_{\mathcal Z} \ar[r]^{b_1}  &  Z \ar[r]^{b_2}    &  W^0\\
 }
\end{xy}
\]

Consider the Cartier divisors $\Delta\subset W'\times W^0$, where $\Delta$ is the restriction of the diagonal of $W'\times W'$, and $D=p\times W^0\subset W'\times W^0$. Set $\phi=\phi_2\circ\phi_1$. By construction,
$\mathcal L=\rho_1^*\mathcal R(-\phi^*\Delta-\phi^*D)$ is a family of theta characteristics on the trivial family $W'\times\Theta_{\mathcal Z}\ra \Theta_{\mathcal Z}$, hence $\mathcal L$ is constant on the connected components of $\Theta_{\mathcal Z}$. Let $q\in W^0$ and $q'\in Z$ such that $b_2(q')=q$. Pick $Z_{q'}=\eta^{-1}(q')$. By Lemma \ref{trivial}, we have that $\pi\col W'\ra Z_{q'}$ is the normalization of $Z_{q'}$ at a tacnode $t_{q'}$ and $\pi^{-1}(t_{q'})=\{p,q\}$. The proof of \cite[Theorem (2.22)]{Ha} and \cite[Section 2e (a)]{Ha} implies that, if $Z$ is a curve, $Z'$ is the normalization at a tacnode of $Z$ with the points $p,q$ over the tacnode and $L$ an odd theta characteristic of $Z'$, then there are exactly 2 even theta characteristics of $Z$ whose pull-back to $Z'$ is $L(p+q)$. Then, if $L_1,\dots,L_N$ are the odd theta characterictics of $W'$ and 
$\eta^{q'}_1,\dots,\eta^{q'}_{2N}\in\Theta_{\mathcal Z}$ are the points representing the even theta characteristics of $Z_{q'}$, then 
by construction $\{\mathcal L|_{\rho^{-1}(\eta^{q'}_1)},\dots,\mathcal L|_{\rho^{-1}(\eta^{q'}_{2N})}\}=\{L_1,L_1,L_2,L_2,\dots,L_N,L_N\}$, where $\rho=\rho_2\circ\rho_1$. Thus $\Theta_{\mathcal Z}$ is the disjoint union of  components $\Theta^1_{\mathcal Z},\dots,\Theta^N_{\mathcal Z}$ such that $\mathcal L|_{\rho^{-1}({\Theta^j_{\mathcal Z}})}$ is the trivial bundle $L_j\ra W'\times\Theta^j_{\mathcal Z}$, for $j=1,\dots, N$. By induction, $h^0(L_j)=1$ and the section of $L_j$ vanishes in smooth points of $W'$, thus we can find $r\in W^0$ such that $h^0(L_j(-r))=1$. Pick $r'\in Z$ such that $b_2(r')=r$. Since $p$ is general, we have $h^0(L_j(p+r))=2+h^0(L_j(-p-r))=2$ and $h^0(L_j(p))=1+h^0(L_j(-p))=1$ and by construction $h^0(L_j(r))=1+h^0(L_j(-r))=2$. In particular, we get a section $s$ of $L_j(p+r)$ vanishing in $p$ and not vanishing in $r$. Thus, $s$ does not descend to a section of the even theta characteristics $R_{2j-1},R_{2j}$ of $Z_{r'}=\eta^{-1}(r')$ whose pull-back to $W'$ is $L_j(p+r)$ and hence $h^0(Z_{r'},R_{2j-1})=h^0(Z_{r'},R_{2j})=0$. Now, $\Theta^j_{\mathcal Z}$ has at most two connected components. If 
$\eta^{r'}_{2j-1},\eta^{r'}_{2j}\in\Theta^j_{\mathcal Z}$ represent $R_{2j-1}, R_{2j}$, the set $\{\eta^{r'}_{2j-1}, \eta^{r'}_{2j}\}$ intersects all the connected components of $\Theta^j_{\mathcal Z}$. By construction, 
$h^0(\mathcal R|_{\rho_2^{-1}(\eta^{r'}_{2j-1})})=h^0(\mathcal R|_{\rho_2^{-1}(\eta^{r'}_{2j})})=0$. Thus, if $q\in W^0$ is general, $b_2(q')=q$ for $q'\in Z$ and $b_1^{-1}(q')\cap\Theta^j_{\mathcal Z}=\{\eta^{q'}_{2j-1},\eta^{q'}_{2j}\}$, by semicontinuity $h^0(Z_{q'}, \mathcal R|_{\rho_2^{-1}(\eta^{q'}_{2j-1})})=h^0(Z_{q'}, \mathcal R|_{\rho_2^{-1}(\eta^{q'}_{2j})})=0$, where $Z_{q'}=\eta^{-1}(q')$ and we are done.

\smallskip

(ii) First of all we show that a hyperplane $H$ of type $0$ with respect to $W$ is a theta hyperplane of $W$ if and only if it cuts $W$ in points with even multiplicity. 
If $H\cdot W=2D_H$, then $\mathcal O_W(D_H)$ is an effective theta characteristic of the canonical curve $W$. From (i), we have that $h^0(\mathcal O_W(D_H))=1$ and Section \ref{sec2.2} implies that $\mathcal O_W(D_H)$ is limit of odd theta characteristics, because the parity of a theta characteristic is preserved in a deformation. Thus $H$ is a theta hyperplane. The converse is obvious. 

We have a set injection:
\[
\SelectTips{cm}{11}
\begin{xy} <16pt,0pt>:
\xymatrix{
\{\text{theta hyperplanes of type 0 of }W\}\UseTips\ar[r] & \{\text{odd theta characteristic of }W\}.
}
\end{xy}
\]

sending $H$ to $\mathcal O_W(D_H)$. If $R$ is an odd theta characteristic of $W,$ let $D$ be the effective divisor of $|R|$ and $H$ be the theta hyperplane cutting $2D$ on $W.$ It follows from (i) that $D$ and hence $H$ do not contain singular points of $W$, thus the injection is also a surjection. 
To prove that $W$ is theta generic, it suffices to show by induction that the set 
$\mathcal H$ of hyperplanes cutting the smooth locus of $W$ in points with even multiplicities is finite. In fact, the theta hyperplanes are contained in $\mathcal H$. 
This is clear if $g=3$ and, by what we have proved, $W$ has a finite number of hyperplanes of type 0 with respect to $W$, cutting $W$ in smooth points with even multiplicity. Now, assume that there is an infinite number of  hyperplanes of $\mathcal H$ containing a fixed set of tacnodes, cusps and nodes. Call $W'$ the curve obtained by projecting $W$ from one of these singular points. The projection of the hyperplanes of $\mathcal H$ gives rise to an infinite set of hyperplanes cutting the smooth locus of $W'$ in points with even multiplicities, which is a contradiction.
\end{proof}

\begin{Lem}\label{cut}
Let $W\subset\Ps^{g-1}$ be an irreducible general canonical curve with tacnodes, cusps and nodes and $H\subset\Ps^{g-1}$ be a hyperplane. Then $H$ is a theta hyperplane of $W$ if and only if $H$ cuts the tacnodes of $W$ with multiplicity $2$ or $4$, the cusps and nodes of $W$ with multiplicity $2$ and smooth points of $W$ with even multiplicities.
\end{Lem}

\begin{proof}
Notice that, if a hyperplane $H$ cuts a tacnode $t\in W$ with multiplicity at least $4$, then $H$ contains the tacnodal tangent $T_t$ of $t$. In fact, let $\pi\col W'\ra W$ be the normalization of $W$ at $t$, with the points $p,q$ lying over $t$. Recall that $\pi^*\omega_W=\omega_{W'}(2p+2q)$.  As in the proof of Lemma \ref{trivial}, if $W'\subset\Ps^g$ is the embedding given by $H^0(\omega_{W'}(2p+2q))$, then $W\subset\Ps^{g-1}$ is the image of $W'$ via the projection 
$\pi_r$ from a point $r$ of the line $\overline{pq}$ with $r\ne p,q$ and a hyperplane $H$ cutting $t$ with multiplicity at least $4$ is the projection via 
$\pi_r$ of a hyperplane of $\Ps^g$ containing either the tangent of $W'$ at $p$ or the tangent of $W'$ at $q$, hence $H$ contains $T_t$.

Let $H$ be as in the statement. If $H$ is of type 0 with respect to $W$, then it is a theta hyperplane by the first part of the proof of Lemma \ref{zero-type} (ii). Let $H$ be of type $(i,j,k,h)$ with respect to $W$, with $t_1,\dots, t_i,c_1,\dots,c_k,n_1,\dots,n_h$ the tacnodes, cusps and nodes contained in $H$. Assume that the tacnodal tangents of $t_1,\dots t_j$ are contained in $H$. 
Let $\pi\col W'\ra W$ be the normalization of $W$ at $t_1,\dots t_j, c_1,\dots c_k,n_1,\dots n_h$ and the partial regularization at $t_{j+1},\dots t_i$, meaning that we get $j$ nodes lying over $t_{j+1},\dots t_i$. If $H$ intersects the smooth locus of $W$ in $2D_H$, we have that $L=\mathcal O_{W'}(\pi^*(D_H))$ is an effective theta characteristic of $W'$, hence $h^0(L)=1$ from Lemma \ref{zero-type} (i). If $W$ has only nodes, consider a general smoothing $\mathcal W\ra B$ of $W$ to theta generic curves and let $B'\ra B$ be a degree 2 covering, totally ramified over $0\in B$. Consider:
\[
\SelectTips{cm}{11}
\begin{xy} <16pt,0pt>:
\xymatrix{
 \mathcal X \UseTips\ar[r]^{\rho}& \mathcal{\widetilde W} \ar[r]& \mathcal W\\}
\end{xy}
\]

where $\mathcal{\widetilde{W}}=\mathcal W\times_B B'$ and $\rho$ is the resolution of the $A_1$ singularities over the nodes $H\cap\text{Sing}(W)$. The central fiber $X$ of $f\col\mathcal X\ra B'$ is a quasistable curve, the union of $W'$ and some exceptional components. Consider a stable spin curve $\xi=(X, G, \alpha)$ of $W$, where $G$ is the gluing of $L\in\Pic(W')$ and the degree 1 bundle over each exceptional component. Thus $h^0(G)=1$, from (\ref{decompos}), hence $G$ is the limit of a family of odd theta characteristics on the smooth fibers of $\mathcal X\ra B'$ and its section is smoothable. Since $\rho$ is an isomorphism away from the exceptional components, by construction the associated family of theta hyperplanes has $H$ as limit hyperplane, hence $H$ is a theta hyperplane of $W$.

Assume that $W$ has either a tacnode or a cusp. Then, Lemma \ref{zero-type} (ii) implies that $W$ is theta-generic. To show that $H$ is a theta hyperplane it suffices to show that it is limit of theta hyperplanes of nodal general curves. Now, consider a one-parameter family of curves $\mathcal W\ra B$ such that the general fiber is a general curve with $i+j+h+k$ nodes and the special fiber is $W$. 
Assume that each one of $t_1,\dots, t_j$ is respectively the limit of 2 nodes and that each one of $t_{j+1},\dots,t_i, c_1,\dots,c_k,n_1,\dots,n_h$ is the limit  of 1 node. The normalization $\mathcal W'$ of the surface $\mathcal W$ gives rise to a smoothing of $W'$. Consider the odd theta characteristic $L$ of $W'$ with $h^0(L)=1$. Of course, $L$ is the limit of a family of odd theta characteristics on the smooth fibers of $\mathcal W'$ and its section $s$ is smoothable. The zero divisor of $s$ is the limit of a family of effective semicanonical divisors on the smooth fibers of $\mathcal W'$, inducing also a family of effective divisors on the nodal fibers of $\mathcal W$. The family of hyperplanes on the nodal fibers of $\mathcal W$, obtained on each fiber as the linear span of the points of the effective divisor and of its $i+j+h+k$ nodes, is a family of theta hyperplanes, because the curve is nodal, and by construction its limit is exactly $H$.  

Conversely, let $H$ be a theta hyperplane of $W$. First of all, $H$ cuts $W$ in points with even multiplicities. If $g=3$, then $H$ cuts the nodes and cusps of the general curve $W$ with multiplicity at most $3$ and a tacnode of $W$ with multiplicity at most $4$, hence we are done. 
Assume that $g\ge 3$. Let $t_1,\dots, t_i,c_1,\dots,c_k,n_1,\dots,n_h$ be the tacnodes, cusps and nodes contained in $H$. Assume that the tacnodal tangents of $t_1,\dots, t_j$ are contained in $H$. In particular, $H$ cuts the remaining tacnodes with multiplicity $2$. Let $\pi\col Z\ra W$ be the normalization of $W$ at $t_1,\dots,t_j,c_1,\dots,c_k,n_1,\dots,n_h$ and the partial regularization at $t_{j+1},\dots,t_i$. 

Let $g_Z$ be the genus of $Z$.  If $g_Z\ge 3$, pick the projection of $W$ from the tacnodal tangents of $t_1,\dots t_j$, the tacnodes $t_{j+1},\dots t_i$ and the cusps and nodes contained in $H$. Then $W$ projects to the canonical model 
$Z\subset\Ps^{g_Z-1}$ of $Z$ and $H$ projects to a hyperplane $H_Z\subset\Ps^{g_Z-1}$ cutting $Z$ only in smooth points. Assume that $H$ cuts $n_s$ with even multiplicity $\alpha_s$, for $s\le h$, $c_s$ with even multiplicity $\beta_s$, for $s\le k$, and $t_s$ with even multiplicity $\gamma_s$, for $s\le j$. In this way, $H_Z\subset\Ps^{g_Z-1}$ cuts the points $u_s,v_s\in Z$ over $n_s$ with multiplicities $\alpha_{s_1},\alpha_{s_2}$ such that $\alpha_{s_1}+\alpha_{s_2}=\alpha_s-2$, the point $w_s\in Z$ over $c_s$ with even multiplicity $\beta_{s_1}=\beta_s-2$ and the points $x_s,y_s\in Z$ over $t_s$, for $s\le j$, with multiplicities $\gamma_{s_1},\gamma_{s_2}$ such that $\gamma_{s_1}+\gamma_{s_2}=\gamma_s-4$. Let $W^{sm}$ be the smooth locus of $W$ and let $D$ be the Cartier divisor such that $2D=H\cdot W^{sm}$. Pick $L=\pi^*\mathcal O_W(D)$. By construction, we have: 
$$L^{\otimes 2}=\omega_Z(-\sum_s (\alpha_{s_1}u_s+\alpha_{s_2}v_s+\beta_{s_1}w_s+\gamma_{s_1}x_s+\gamma_{s_2}y_s)).$$ In particular, if we set:
$$M=L([\frac{\alpha_{s_1}+1}{2}]u_s+[\frac{\alpha_{s_2}+1}{2}]v_s+\frac{\beta_{s_1}}{2}w_s+[\frac{\gamma_{s_1}+1}{2}]x_s+[\frac{\gamma_{s_2}+1}{2}]y_s),$$ then $M^{\otimes 2}=\omega_Z(\sum_{s\in I}(u_s+v_s)+\sum_{s\in J}(x_s+y_s))$, where $I$ is the set of the indices such that $\alpha_{s_1}$ is odd and $J$ is the set of the indices such that $\gamma_{s_1}$ is odd. Assume that either $I\ne\emptyset$ or $J\ne\emptyset$ and let $Z'$ be the general curve obtained by the nodal identifications of the pairs of points $u_s,v_s\in Z$ indexed by $I$ and the pairs of points $x_s,y_s\in Z$ indexed by $J$. We get a theta characteristic $L'$ on $Z'$ whose pull-back to $Z$ is $M$. By construction, a section of $M$ vanishes on $u_s,v_s$, if $s\in I$ and on $x_s,y_s$, if $s\in J$, hence $L'$ has a section vanishing on some node, which contradicts Lemma \ref{zero-type} (i). Thus $I=J=\emptyset$ and $M^{\otimes 2}=\omega_Z$. Notice that $h^0(M)\ge 1$, hence by Lemma \ref{zero-type} (i) we have $h^0(M)=1$. The section of $M$ does not vanish on $u_s,v_s,w_s,x_s,y_s$, because they are general points of $Z$, hence $\alpha_{s_1}=\alpha_{s_2}=\beta_{s_1}=\gamma_{s_1}=\gamma_{s_1}=0$. In this way, we have $\alpha_s=2$, for $s\le h$, $\beta_s=2$, for $s\le k$, and $\gamma_s=4$, for $s\le j$.

Assume that $g_Z\le 2$. If $H$ cuts a node $n_1$ with multiplicity $m\ge 4$, there exists a projection from a suitable subset of the set of tacnodal tangents of $t_1\dots,t_j$ and the points $t_{j+1},\dots,t_i, c_1,\dots c_k, n_2,\dots n_h$, such that $W$ projects to a general plane quartic with a node $n_1$. Now, $H$ projects to a line intersecting $n_1$ with multiplicity $m$, which is a contradiction. A similar argument works if $H$ contains a cusp with multiplicity $m\ge 4$ or a tacnode with multiplicity $m>4$.
\end{proof}

Recall the definition of $N_g:=2^{g-1}(2^g-1)$ and $N^+_g:=2^{g-1}(2^g+1)$.

\begin{Thm}\label{hyper}
Let $W$ be an irreducible general canonical curve with $\tau$ tacnodes, $\gamma$ cusps and $\delta$ nodes of genus $g\ge 3.$ Let $\widetilde g=g-\delta-\gamma-2\tau$ be the genus of the normalization of $W.$ 

If $j<i$ or $h\ne\delta,$ then: $$t^j_{ikh}(W)=2^{2\widetilde{g}+\tau-j+\delta-h-1}\binom{\tau}{i}\binom{i}{j}\binom{\gamma}{k}\binom{\delta}{h}.$$

If $i=j$ and $h=\delta,$ then: $$t^i_{ik\delta}(W)=\left\{\begin{array}{ll}\displaystyle   2^{\tau-i}\binom{\tau}{i}\binom{\gamma}{k} N_{\widetilde{g}}& \text{if \ }\tau -i+ \gamma -k\equiv 0 \;(2) \\ \\ 
\displaystyle 2^{\tau-i}\binom{\tau}{i}\binom{\gamma}{k} N^+_{\widetilde{g}}&\text{if \ }\tau-i+\gamma-k\equiv 1 \; (2);.\end{array}\right.$$
\end{Thm}

\begin{proof}
The proof is by induction on $g.$ The formulas hold if $W$ is smooth and if $g=3$, by \cite[3.2]{CapSerbi}. 
Consider the case $(i,j,k,h)\ne (0,0,0,0).$ Call $W'$ the projection of $W$ form a singular point $s$: since $g\ge 4$ we can project at least one time. Lemma \ref{cut} implies that $H\ni s$ is a theta hyperplane of $W$ if and only if the projection of $H$ is a theta hyperplane of $W'$. Denote by $W^g_{\tau\gamma\delta}$ a genus $g$ irreducible general canonical curve with $\tau$ tacnodes, $\gamma$ cusps and $\delta$ nodes. Thus $t^j_{ikh}(W^g_{\tau\gamma\delta})$ is obtained by multiplying the cardinality of the set $\mathcal T$ of theta hyperplanes of $W^g_{\tau\gamma\delta}$ of type $(i,j,k,h)$ containing a fixed set of $i$ tacnodes, $j$ tacnodal tangents, $k$ cusps and $h$ nodes and the number $\binom{\tau}{i}\binom{i}{j}\binom{\gamma}{k}\binom{\delta}{h}$ of all possible sets of such singularities.

We compute the cardinality of $\mathcal T$. If $j<i,$ we project the curve from a tacnode $t$ contained in the theta hyperplanes of $\mathcal T$ and whose tacnodal tangent is not contained in the hyperplanes. The projected curve $W^{g-1}_{\tau-1,\gamma,\delta+1}$ has genus $g-1$ and we can apply the induction. By the first argument of the proof of Lemma \ref{cut}, the theta hyperplanes of $\mathcal T$ intersects $t$ with multiplicity $2$, hence they project to the theta hyperplanes of $W^{g-1}_{\tau-1,\gamma,\delta+1}$ of type $(i-1,j,k,h)$ containing a fixed set of $i-1$ tacnodes, $j$ tacnodal tangents, $k$ cusps and $h$ nodes. Then:
$$|\mathcal T|=\frac{t^j_{ikh}(W^g_{\tau\gamma\delta})}{\binom{\tau}{i}\binom{i}{j}\binom{\gamma}{k}\binom{\delta}{h}}=\frac{t^j_{i-1,k,h}(W^{g-1}_{\tau-1,\gamma,\delta+1})}{\binom{\tau-1}{i-1}\binom{i-1}{j}\binom{\gamma}{k}\binom{\delta+1}{h}}.$$ 

The curves $W^{g-1}_{\tau-1,\gamma,\delta+1}$ and $W^g_{\tau\gamma\delta}$ have the same normalization and $\delta+1\ne h,$ thus we are done because by induction: 
$$t^j_{i-1,k,h}(W^{g-1}_{\tau-1,\gamma,\delta+1})=\binom{\tau-1}{i-1}\binom{i-1}{j}\binom{\gamma}{k}\binom{\delta+1}{h}2^{2\widetilde{g}+\tau-j+\delta-h-1}.$$

If $i=j\ne 0$ and $\delta=h,$ we project the curve from a tacnode $t$ contained in the theta hyperplanes of $\mathcal T$. The tacnode projects to a node $n$. Since the hyperplanes of $\mathcal T$ contain the tacnodal tangent of $t$, they project to hyperplanes containing $n$. Thus we have: 
$$|\mathcal T|=\frac{t^i_{ikh}(W^g_{\tau\gamma\delta})}{\binom{\tau}{i}\binom{\gamma}{k}}=\frac{t^{i-1}_{i-1,k,h+1}(W^{g-1}_{\tau-1,\gamma,\delta+1})}{\binom{\tau-1}{i-1}\binom{\gamma}{k}}.$$ 

Since $\delta+1=h+1,$ we are done because by induction we have: 
$$t^{i-1}_{i-1,k,h+1}(W^{g-1}_{\tau-1,\gamma,\delta+1})=\binom{\tau-1}{i-1}\binom{\gamma}{k}2^{\tau-1-i+1}N_{\widetilde{g}}$$
if $\tau-i+\gamma-k\equiv 0 \; (2)$ 
and $$t^{i-1}_{i-1,k,h+1}(W^{g-1}_{\tau-1,\gamma,\delta+1})=\binom{\tau-1}{i-1}\binom{\gamma}{k}2^{\tau-1-i+1}N^{+}_{\widetilde{g}}$$ 
if $\tau-i+\gamma-k\equiv 1 \; (2).$

In the other cases, we argue as before by projecting the curve from a tacnode if $i\ne 0,$ from a node if $i=0$ and $h\ne 0,$ from a cusp if $i=h=0.$ 

By Lemma \ref{zero-type} (ii), the number $t^0_{000}(W)$ is the number of odd theta characteristics of $W$, hence it is given by \cite[Corollary 2.7, Corollary 2.8]{Ha}, where $k=\tau+\delta,$ as explained in \cite[Section 2b]{Ha}. If $\delta\ne 0$, it follows from \cite[Theorem 2.12]{Ha} that $t^0_{000}(W)$ is given by \cite[Corollary 2.7]{Ha} and we are done. 
If $\delta=0$, again from \cite[Theorem 2.12]{Ha} we have that $t^0_{000}(W)$ is given by \cite[Corollary 2.8]{Ha}. 
As explained in \cite[Section 2e and Theorem 2.22]{Ha}, the formula (2.9) (respectively (2.10)) of \cite[Corollary 2.8]{Ha} holds if $\tau+\gamma$ is even (respectively odd). 
\end{proof}

\section{The multiplicity of a theta hyperplane}\label{sec5}

We complete the description of the zero-dimensional scheme of theta hyperplanes of irreducible curves, computing the multiplicities of its points. We will denote by $\omega_f$ the relative dualizing sheaf of a family of curves $f\col\mathcal W\ra B.$

\begin{Lem}\label{redcal}
Let $W$ be an irreducible curve of genus at least 2 with normalization $W^\nu.$ Let $\mathcal W\ra B$ be a general smoothing of $W.$ Let $C$ be the central fiber of the stable reduction $\mathcal C$ of $\mathcal W.$
\itemize
\item[(i)] 
Assume that $W$ has exactly $\gamma$ cusps as singularities. Let $b\col B'\ra B$ be the base change of order $6$ totally ramified over $0\in B.$ Then $\mathcal C$ is a smooth surface and there are smooth elliptic curves $F_1,\dots,F_\gamma,$ so that the dual graph of $C$ is given by:
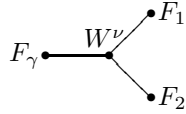
\[
\begin{xy} <16pt,0pt>:
(0,0)*{\scriptstyle\bullet}="a"; 
(1,-1)*{\scriptstyle\bullet}="b";
(1,1)*{\scriptstyle\bullet}="d";
(-1.5,0)*{\scriptstyle\bullet}="f";
"a"+0;"b"+0**\dir{-};
"a"+0;"d"+0**\dir{-};
"a"+0;"f"+0**\dir{-};
"b"+(0.5,0)*{F_2};
"d"+(0.5,0)*{F_1}; 
"f"+(-0.5,0)*{F_\gamma}; 
"a"+(-0.1,0.4)*{W^\nu};
\end{xy}
\]

Assume that $W$ has exactly $\tau$ tacnodes as singularities. Let $b\col B'\ra B$ be the base change of order $4$ totally ramified over $0\in B.$ Then $\mathcal C$ is a smooth surface and there are smooth elliptic curves $F_1,\dots,F_\tau,$ so that the dual graph of $C$ is given by:
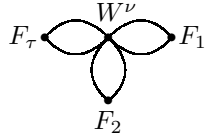
\[
\begin{xy} <16pt,0pt>:
(0,0)*{\scriptstyle\bullet}="a"; 
(0,-1.5)*{\scriptstyle\bullet}="b";
(1.5,0)*{\scriptstyle\bullet}="d";
(-1.5,0)*{\scriptstyle\bullet}="f";
"a"+0;"b"+0**\crv{(-0.8,-0.8)};
"a"+0;"b"+0**\crv{(0.8,-0.8)}; 
"a"+0;"d"+0**\crv{(0.8,-0.8)};
"a"+0;"d"+0**\crv{(0.8,0.8)}; 
"a"+0;"f"+0**\crv{(-0.8,-0.8)};
"a"+0;"f"+0**\crv{(-0.8,0.8)}; 
"b"+(0,-0.5)*{F_2};
"d"+(0.5,0)*{F_1}; 
"f"+(-0.5,0)*{F_\tau}; 
"a"+(0.2,0,6)*{W^\nu};
\end{xy}
\]
\item[(ii)]
Let $\mathcal W$ be as in $(\first),$ with stable reduction $f\col\mathcal C\ra B'.$  Consider the Cartier divisor of $\mathcal C$ given by $F:=\sum F_h,$ the sum of all the elliptic components. Consider $h\col\mathcal W'=\mathcal W\times_{B} B'\ra B'.$ Then $\mathcal C$ is endowed with a $B'$-morphism $\phi\col\mathcal C\ra\mathcal W'$ such that:
\begin{equation}\label{PB}
\phi^*(\omega_h)\simeq\omega_f(F).
\end{equation}
\end{Lem}

\begin{proof}
(i) We argue as in \cite[Theorem III-10.1]{BPV} and \cite[Example pag.122]{HaMo}. This proof works, up to replace $B$ by some open subset containing $0$. Since $\mathcal W$ is general, it is a smooth surface. Let $\mathcal{\overline{W}}$ be obtained by blowing-up $\mathcal W$ three times in correspondence of each cusp, so that the reduced special fiber has normal crossings. Take a base change $b_1\col B_1\ra B$ of order $2$ totally ramified over $0\in B$ and the normalization $\mathcal W_1$ of $\mathcal{\overline{W}}\times_{B_1} B.$ As explained in \cite{HaMo}, $\mathcal W_1$ is the double cover of $\mathcal W$ branched along the irreducible components of the special fiber of $\mathcal{\overline W},$ appearing with odd multiplicities. Then $\mathcal W_1$ is a smooth surface, because the branch divisor is smooth. Take the base change $b_2\col B'\ra B_1$ of order $3$ totally ramified over $0\in B_1$ and the normalization $\mathcal C'$ of $\mathcal W_1\times_{B_1}B'.$ Then $\mathcal C'$ is the triple cover of $\mathcal W_1$ ramified along the irreducible components of the special fiber, appearing with multiplicities not divisible by $3.$ Then $\mathcal C'$ is a smooth surface because the branch divisor is smooth. The components of the special fiber of $\mathcal C'$ are $W^\nu,$ $\gamma$ elliptic curves intersecting transversally $W^\nu$ in one point, and $(-1)$-curves. Then $\mathcal C$ is obtained by contracting the $(-1)$-curves contained in the special fiber of $\mathcal C'\ra B'.$

The tacnodal case is similar, combining two base changes of order $2$ totally ramified over $0.$ 

(ii) Let $\mathcal C'$ be as in (i). By the universal property of the fiber products, we have a $B'$-morphism from $\mathcal C'$ to $\mathcal W',$ factorizing through the $B'$-relative minimal model $\mathcal C$ of $\mathcal C'.$ 
We get the diagram:
\[
\SelectTips{cm}{11}
\begin{xy} <16pt,0pt>:
\xymatrix{
 \mathcal C \ar[dr]_{f}\UseTips\ar[r]^{\phi}& \mathcal W' \ar[r]\ar[d]^{h} & \mathcal W\ar[d]\\
     &  B'\ar[r]^{b}     &  B\\}
\end{xy}
\]

Since $\phi$ is an isomorphism away from the special fibers, it follows that $\omega_f$ and $\phi^*(\omega_h)$ differ by a divisor $D$ of $\mathcal C$ supported on components of $C.$ If $\nu\col W^\nu\ra W$ is the normalization, it follows that $\phi^*(\omega_h)\otimes\mathcal O_{W^\nu}\simeq\nu^*(\omega_W)\simeq\omega_{W^\nu}(2\sum(F_h\cap F_h^c))$ and hence $D\sim F.$ Thus we get the relation (\ref{PB}).
\end{proof}

\begin{Def}
We call the elliptic curves $F_h,$ appearing in the stable reduction, \emph{elliptic tails}.
\end{Def}

\begin{Rem}\label{sing-des}
Consider the stable reduction of a general smoothing of a tacnodal curve, as in Lemma \ref{redcal}. It is easy to see that an elliptic tail $F_h$ is a double cover $\psi\col F_h\ra\Ps^1,$ branched at $0,1,\infty,-1,$ with the points $F_h\cap F^c_h$ lying over $0,\infty.$ It is the elliptic curve with $j$-invariant $j=1728$.
\end{Rem}

\begin{Lem}\label{aut-exc}
Let $\mathcal C$ be the stable reduction of a general smoothing of an irreducible curve of genus at least 2 with a tacnode. Let $F$ be the elliptic tail of $\mathcal C$ over the tacnode and set $F\cap F^c=\{p,q\}.$ There exists a group $\mathcal A=\{id,\gamma_1,\gamma_2,\gamma_3\}$ of automorphisms of $F,$ fixing $p$ and $q$ and such that, if $G_1$ and $G_2$ are square roots of $\mathcal O_F(p+q),$ then
$\gamma^*G_1\simeq G_2$ for some $\gamma$ in $\mathcal A.$
\end{Lem}

\begin{proof} 
By Remark \ref{sing-des}, $F$ is a double cover $\psi\col F\ra\Ps^1$ branched over $0,1,\infty,-1$ and $p,q$ lie over $0,\infty.$ Thus $F$ admits an involution $\gamma_1,$ exchanging the ramification points $r_1,r_2$ over $1,-1$ and whose set of fixed points is $\{p,q\}.$ Let $\gamma_2$ be the involution of $F$ associated to $\psi$ and $\gamma_3:=\gamma_2\circ\gamma_1.$ The four square roots of $\mathcal O_F(p+q)$ are effective. Pick the distinct effective divisors $D_i,$ for $i=1,\dots,4,$ such that $2D_i\sim p+q.$ Of course $D_i\ne p,q,r_1,r_2.$ Since $2\gamma_j^*(D_i)=\gamma_j^*(2 D_i)\sim \gamma_j^*(p+q)=p+q,$ for $j=1,2,$ we get, up to the order, $D_1=\gamma_1^*(D_3), D_2=\gamma_1^*(D_4)$ and
 $D_2=\gamma_2^*(D_1), D_3=\gamma_2^*(D_4).$ 
 \end{proof}

\begin{Lem}\label{secvan}
Let $C$ be the nodal union of a canonical irreducible general curve $W$ with $\tau$ tacnodes, $\gamma$ cusps and $\delta$ nodes and $m$ smooth elliptic curves $F_1,\dots,F_m$ with $F_h\cap F_k=\emptyset$ for $h\ne k$. Assume that $F_h\cap F_h^c$ are two general points of $W$. Let $R$ be an odd theta characteristic of $C$. Then $h^0(R)=1$ and the non-trivial section of $R$ vanishes on smooth points of $C$.
\end{Lem}

\begin{proof}
The proof is by induction on $m$. The case $m=0$ is Lemma \ref{zero-type}. 
Let $g_W$ be the genus of the normalization of $W$. Set $F_h\cap F_h^c=\{p_h, q_h\}$. Consider the effective degree 1 Cartier divisors 
$D_{h1},D_{h2},D_{h3},D_{h4}$ of $F_h$ whose double is in $|p_h+q_h|$. Abusing notation, we see $D_{hj}$ as an effective degree 1 Cartier divisor of $C$ supported in a smooth point of $C$ contained in $F_h$. Let $W'$ be the irreducible general curve obtained by the nodal identification of the two points of $F_h\cap F^c_h$ for every $h$. Since $W'$ has at least one node, it follows from \cite[Corollary (2.7), Theorem (2.12)]{Ha} that $W'$ has $N_W=2^{2g_W+k-1}$  odd theta characteristics, where $k=m+\delta+\tau$. By induction, they have one section whose divisor of zeroes is supported in smooth points of $W'$. Since $W$ is the normalization of $W'$ at the $m$ new nodes, by pull-back we get $N_W$ effective Cartier divisors of $W$ whose double is in $|\omega_W(\sum_{1\le h\le m}(p_h+q_h))|$, which we can see as effective Cartier divisors of $C$ supported on smooth points of $C$ contained in $W$. Consider the sum of each one of these $N_W$ effective Cartier divisors with one $D_{hj}$ for every elliptic curve $F_h$. We get a set of $4^m N_W$ effective Cartier divisors of $C$, supported on smooth points of $C$. By construction, the set $\mathcal S$ of the associated line bundles is a set of effective theta characteristics of $C$. Now, $C$ has at least one node and its normalization has genus $g_W+m$, thus \cite[Corollary (2.7)]{Ha} implies that $C$ has exactly $2^{2(g_W+m)+k-1}=4^m N_W$ odd theta characteristics. If we show that $h^0(R)=1$ for every $R\in\mathcal S$, then 
$\mathcal S$ is the set of the $4^m N_W$ odd theta characteristics of $C$ and we are done.

Let $R\in\mathcal S$. We have $H^0(R)\subseteq H^0(R|_{F_1})\oplus H^0(R|_{F_1^c})$ and $h^0(R|_{F_1})=1$. If $h^0(R)\ge 2$, then there exists $0\ne s\in H^0(R|_{F_1^c})$ such that $(0, s)$ descends to a section of $R$. In particular the section $s$ vanishes on $F_1\cap F_1^c$. Consider the curve $W_1$ obtained from $F_1^c$ by the nodal identification of $F_1\cap F_1^c$. Call $n$ the new node. From \cite[Theorem (2.14)]{Ha} we get exactly one odd theta characteristic $R_1$ of $W_1$ whose pull-back to $F_1^c$ is $R|_{F_1^c}$ and $s$ descends to a section of $R_1$ vanishing on $n$, which is not possible by induction.
\end{proof}

\subsection{Curves of twisted spin curves} Let $\mathcal W\ra B$ be a general smoothing of a curve $W$ as in Lemma \ref{redcal}. Pick its stable reduction $f\col\mathcal C\ra B'.$ 
For a divisor $D$ of $\mathcal C,$ supported on $C$, set $\mathcal N_D:=\omega_f(D).$ Pick the moduli space:
\[
\SelectTips{cm}{11}
\begin{xy} <16pt,0pt>:
\xymatrix{
\overline{S}_f(\mathcal N_D)\UseTips\ar[r] & B'
}
\end{xy}
\]

of \cite[Theorem 2.4.1.]{CapCasCorn}. Let $S^-_{\omega_f^*}$ be the open subscheme of $\overline{S}_f(\mathcal N_D)$ of odd theta characteristics of smooth fibers of $\mathcal C\ra B'.$ Let $S^-_{\mathcal N_D}$ be the closure of $S^-_{\omega_f^*}$ in $\overline{S}_f(\mathcal N_D).$ The curves $S^-_{\mathcal N_D}$ are  birational, as $D$ varies. Thus they have the same curve as normalization, which we denote by:
\[
\SelectTips{cm}{11}
\begin{xy} <16pt,0pt>:
\xymatrix{
\nu_D\col S^\nu_f \UseTips\ar[r] & S^-_{\mathcal N_D}.
}
\end{xy}
\]

Let $J_{\mathcal W'}$ be the curve of theta hyperplanes of Section \ref{sec4}. We get a rational map:
\[
\SelectTips{cm}{11}
\begin{xy} <16pt,0pt>:
\xymatrix{
\mu_D\col S^-_{\mathcal N_D} \UseTips\ar@{.>}[r] & J_{\mathcal W'}
}
\end{xy}
\]

which is an isomorphism away from the central fiber. Since $S^\nu_f$ is smooth, we get a morphism:
\begin{equation}\label{nonst-spin}
\SelectTips{cm}{11}
\begin{xy} <16pt,0pt>:
\xymatrix{
\psi\col S^\nu_f \UseTips\ar[r] & J_{\mathcal W'}.
}
\end{xy}
\end{equation}

The morphism $\psi$ generically associates to an odd theta characteristic of a smooth curve, the theta hyperplane induced by its unique section. Over the special fiber, the morphism will be described by looking at the behaviour of the smoothable sections of twisted spin curves.
 
\begin{Thm}\label{mult}
Let $W$ be an irreducible general canonical curve of genus $g$ whose singular points are tacnodes and cusps. Then the multiplicity of a theta hyperplane of type $(i, j, k)$ is $4^{i-j} 6^j 3^k.$
\end{Thm}

\begin{proof}
Assume that the singular points of $W$ are exactly $\tau$ tacnodes $t_1,\dots,t_\tau.$ By Theorem \ref{hyper}, $W$ is theta-generic. Let $\mathcal W\ra B$ be a general smoothing of $W$ to theta-generic curves and $f\col\mathcal C\ra B'$ be its stable reduction as in Lemma \ref{redcal} with central fiber $C.$ Set $\{n_{h1},n_{h2}\}:=F_h\cap F^c_h.$

\[
\begin{xy} <16pt,0pt>:
(0,0)*{\scriptstyle}="a"; 
"a"+(-10,3);"a"+(-2,3)**\crv{(-6,0.5)};
"a"+(-10,2);"a"+(-8.5,1.5)**\crv{(-8,4.5)}; 
"a"+(-2,2);"a"+(-3.5,1.5)**\crv{(-4,4.5)};
"a"+(-6.8,1);"a"+(-5.2,1)**\crv{(-6,4)}; 
"a"+(1,3);"a"+(2,2)**\crv{(2,3)}; 
"a"+(2,2);"a"+(3,3)**\crv{(2,3)}; 
"a"+(2,2);"a"+(2,1.8)**\crv{(2,2)}; 
"a"+(2,1.8);"a"+(2,1)**\crv{(1.2,1)}; 
"a"+(2,1.8);"a"+(2,1)**\crv{(2.8,1)}; 
"a"+(3,3);"a"+(4,2)**\crv{(4,3)}; 
"a"+(4,2);"a"+(5,3)**\crv{(4,3)}; 
"a"+(4,2);"a"+(4,1.8)**\crv{(4,2)}; 
"a"+(4,1.8);"a"+(4,1)**\crv{(3.2,1)}; 
"a"+(4,1.8);"a"+(4,1)**\crv{(4.8,1)}; 
"a"+(5,3);"a"+(6,2)**\crv{(6,3)}; 
"a"+(6,2);"a"+(7,3)**\crv{(6,3)}; 
"a"+(6,2);"a"+(6,1.8)**\crv{(6,2)}; 
"a"+(6,1.8);"a"+(6,1)**\crv{(5.2,1)}; 
"a"+(6,1.8);"a"+(6,1)**\crv{(6.8,1)}; 
"a"+(-6,3.5)*{C};
"a"+(4.1,3.5)*{W};
"a"+(-10.2,3.2)*{W^\nu};
"a"+(-10.3,2)*{F_1}; 
"a"+(-6,-0.3)*{\text{The curve $C$ for $\tau=3$}}; 
"a"+(4,-0.3)*{\text{The curve $W$ for $\tau=3$}}; 
"a"+(-7.1,0.9)*{F_2};
"a"+(-3.8,1.4)*{F_3};
"a"+(5.7,2.2)*{t_3};
"a"+(3.7,2.2)*{t_2};
"a"+(1.7,2.2)*{t_1};
\end{xy}
\]

For every divisor $D$ of $\mathcal C,$ consider the diagram:
\[
\SelectTips{cm}{11}
\begin{xy} <16pt,0pt>:
\xymatrix{
 S^\nu_f \ar[d]_{\nu_0}\UseTips\ar[r]^{\nu_D}\ar[dr]^{\psi}& S^-_{\mathcal N_D} \ar@{.>}[d]^{\mu_D}\\
S^-_{\omega_f} \ar@{..>}[r]^{\mu_0} & J_{\mathcal W'}
}
\end{xy}
\]

such that $\mu_D\circ\nu_D=\psi,$ where $\mu_D$ is defined. 
We compute the multiplicities of the special fiber of $J_{\mathcal W'}\ra B'$ by describing the map $\psi.$

\smallskip

\noindent
\emph{First Step: the reduction to twisted spin curves}. For every subset $H\subseteq\{1,\dots,\tau\}$, let $\mathcal S_H$ be the set of stable spin curves in $S^-_{\omega_f}$ supported on the curve $X$ obtained by blowing-up in $C$ the nodes $n_{h1},n_{h2}$ for $h\in H$. Set $\mathcal S^\nu_H=\nu_0^{-1}(\mathcal S_H)$ and $D_H=\sum_{h\in H}F_h$. The goal of the first step is:
\begin{itemize}
\item[(1)]
to describe explicitly $\nu_{D_H}(\nu_0^{-1}(\xi))$ for every $\xi\in\mathcal S_H$,
\item[(2)]
to show that $\nu_{D_H}$ is an isomorphism over $S^\nu_H$.
\end{itemize}

Pick $\xi\in S^-_{\omega_f},$ supported on a blow-up $X$ of $C.$ Let $(X,G,\alpha)$ be a representative of the isomorphism class of $\xi.$ The possible blow-ups of $C$ are described in the Example \ref{banana}. Assume that the nodes which are blown-up to get $X$ are $\{n_{h1},n_{h2}\}$ for $h=1,\dots,j$ and $h=i+1,\dots,\tau$ where $0\le j\le i\le\tau.$ Let $E_{h1}, E_{h2}$ be the exceptional components of $X$ connecting $F_h$ to $W^\nu$. If $A_X$ is the graph associated to $X,$ as in Definition \ref{admgraph}, then $A_X=\Sigma_X.$ The dots of the drawing mean that there are loops from $F_1$ to $F_j$ and from $F_{i+1}$ to $F_\tau$.
\[
\begin{xy} <16pt,0pt>:
(0,0)*{\scriptstyle\bullet}="a"; 
(0,1.5)*{\scriptstyle\bullet}="c";
(0,-1.5)*{\scriptstyle\bullet}="b";
(1.5,0)*{\scriptstyle\bullet}="d";
(-1.5,0)*{\scriptstyle\bullet}="f";
"a"+(0,-3)*{\text{The graph }A_X};
"a"+0;"c"+0**\crv{(-0.8,0.8)};
"a"+0;"c"+0**\crv{(0.8,0.8)}; 
"a"+0;"b"+0**\crv{(-0.8,-0.8)};
"a"+0;"b"+0**\crv{(0.8,-0.8)}; 
"a"+0;"d"+0**\crv{(0.8,-0.8)};
"a"+0;"d"+0**\crv{(0.8,0.8)}; 
"a"+0;"f"+0**\crv{(-0.8,-0.8)};
"a"+0;"f"+0**\crv{(-0.8,0.8)}; 
"b"+(0,-0.5)*{F_j};
"d"+(0.5,0)*{F_1}; 
"c"+(0,0.5)*{F_\tau};
"a"+(-1,1)*{\dots};
"a"+(1,-1)*{\dots};
"f"+(-0.7,0)*{F_{i+1}}; 
"a"+(0.9,0.7)*{\scriptstyle{W^\nu}};
\end{xy}
\begin{xy} <16pt,0pt>:
(0,0)*{\scriptstyle}="a";
"a"+(6,1.3)*{X}; 
"a"+(6,-3)*{\text{The curve $X$ for $\tau=3$ and $i=2j=2$}}; 
"a"+(2,1);"a"+(10,1)**\crv{(6,-1.5)};
"a"+(5.2,-1);"a"+(6.8,-1)**\crv{(6,2)}; 
"a"+(3,1);"a"+(1.5,-0.5)**\dir{-};
"a"+(4,0.5);"a"+(3,-1)**\dir{-};
"a"+(1,0.5);"a"+(4,-1)**\crv{(2.5,0)};
"a"+(9,1);"a"+(10.5,-0.5)**\dir{-};
"a"+(8,0.5);"a"+(9,-1)**\dir{-};
"a"+(8,-1);"a"+(11,0.5)**\crv{(9.5,0)};
"a"+(5,-1.2)*{F_2}; 
"a"+(1.9,1.2)*{W^\nu}; 
"a"+(0.7,0.5)*{F_1}; 
"a"+(1.3,-0.7)*{E_{11}}; 
"a"+(2.8,-1.2)*{E_{12}}; 
"a"+(11.4,0.5)*{F_3}; 
"a"+(11,-0.7)*{E_{32}}; 
"a"+(9.5,-1.2)*{E_{31}}; 
\end{xy}
\]

If no otherwise specified, in the First and Second Step $\xi$ will be fixed. 

Assume that:
\begin{equation}\label{restr}
G|_{F_h}\text{ is even for }1\le h\le j;\,\,G|_{F_h}=\mathcal O_{F_h} \text{ for }i<h\le\tau.
\end{equation}

Of course, $G|_{W^\nu\cup F_{j+1}\dots\cup F_i}$ is a theta characteristic of $W^\nu\cup F_{j+1}\dots\cup F_i$.

To describe the map $\psi$, we choose other representatives in the equivalence class of $\xi$ as follows. By Lemma \ref{redcal}, we have that $\mathcal C$ is smooth. Define the divisor of $\mathcal C:$ 
\begin{equation}\label{Cart-def}
D:=\underset{1\le h\le  j}{\sum} F_h+\underset{i<h\le\tau}{\sum}F_h.
\end{equation}

Pick a $D$-twisted spin curve $(C,L),$ where $L\in\Pic C$ satisfies $L^{\otimes 2}\simeq\omega_C\otimes\mathcal O_f(D)$ and has restrictions:
\begin{equation}\label{restr'}
L|_{F_h}=G|_{F_h} \text{ for }1\le h\le j;\,\,\,\,L|_{F_h}=G|_{F_h}=\mathcal O_{F_h} \text{ for }i<h\le\tau
\end{equation}
\begin{equation*}
\textstyle
L|_{W^\nu\cup F_{j+1}\cup\dots F_i}=G|_{W^\nu\cup F_{j+1}\dots\cup F_i}(\sum_{1\le h\le j}(n_{h1}+n_{h2})+\sum_{i<h\le\tau}(n_{h1}+n_{h2}))
\end{equation*}
 
Since $b_1(A_X)=\tau-i+j,$ there are $2^{\tau-i+j}$ possible gluings, giving rise to a set $\mathcal S_\xi$ of $2^{\tau-i+j}$ line bundles $L$ as above. 

We claim that for every $L\in\mathcal S_\xi$, there exists a representative $(X, G, \alpha)$ of $\xi$ such that $L$ and $G$ are limits of the same family of theta characteristics. In fact, consider: 
\[
\SelectTips{cm}{11}
\begin{xy} <16pt,0pt>:
\xymatrix{
 \mathcal X \UseTips\ar[r]^{\rho_2}& \mathcal{\widetilde C} \ar[r]^{\rho_1}& \mathcal C\\}
\end{xy}
\]

where $\rho_1$ is a double cover ramified over the central fiber, and $\rho_2$ is the resolution of the $A_1$ singularities over $n_{h1}, n_{h2}$ for $h\le j$ and $h>i$. Notice that $\mathcal X$ is a smoothing of $X$. Set $\rho=\rho_2\circ\rho_1$. Now, by Section \ref{sec2.2}, we have that $L$ extends locally (analitically) to a line bundle 
$\mathcal L$ on $\mathcal C$, which is a family of theta characteristics away from the special fiber. Consider:
\begin{equation*}
\textstyle
\mathcal G:=\rho^*\mathcal L(\sum_{h\le j,h>i}(-F_h-E_{h1}-E_{h2})).
\end{equation*}

By construction, the restriction of $\mathcal G$ to $X$ is a line bundle appearing in a representative of $\xi$ and the claim follows.

Now, $S^\nu_f$ is complete over $B',$ then $\emptyset\ne\nu^{-1}_D(C, L)\cap \nu^{-1}_0(\xi)=(C, L)\in S^\nu_f,$ because, by the construction of $S^-_{\mathcal N_D}$ given in \cite{CapCasCorn}, $S^-_{\mathcal N_D}$ is smooth in $(C, L).$ Since $\Sigma_X=A_X,$  Example \ref{banana} implies that the multiplicity of $\xi$ in $S^-_C$ is $2^{\tau-i+j}$. There are $N_g$ odd spin curves of $C,$ with multiplicity, then $\{\nu_D^{-1}(C,L):\forall\,L\in\mathcal S_\xi,\forall\,\xi\in S^-_C\}$ is a subset of cardinality $N_g$ of the fiber of $S^\nu_f$ over $0\in B',$ hence it is the entire fiber. The conclusion is that $\nu_D(\nu_0^{-1}(\xi))=\{(C,L) \;|\; L\in\mathcal S_\xi\}$ and $\nu_D$ is an isomorphism over $\nu_0^{-1}(\xi)$. 

Now, $S^-_{\mathcal N_D}$ is smooth at $(C,L),$ hence $\mu_D$ is defined on $(C,L).$ 
Then we can write:  $$\{\psi(\xi') \;|\;\xi'\in\nu_0^{-1}(\xi)\}=\{\mu_D(C,L) \;|\; L\in\mathcal S_\xi\}.$$

\noindent
\emph{Second Step: the smoothable sections of the line bundles of $\mathcal S_\xi$}. Let $(C,L),$ with $L\in\mathcal S_\xi.$ Pick a representative $\xi=(X,G,\alpha),$ such that $L$ and $G$ are limit of the same family of theta characteristics. Let $\mathcal L$ be the line bundle on $\mathcal C$ extending $L$ to a family of theta characteristics away from the special fiber. Notice that $f\col\mathcal C\ra B'$ is a smoothing to theta-generic curves and $L$ and $G$ are limit of the same family of theta characteristics, then there is a unique $\mathcal L$-smoothable section of $L.$ 

We describe the behavior of this section. As in Lemma \ref{redcal}, we can consider the stable reduction $g\col\mathcal Y\ra B'$ of $\mathcal W$ at $t_1,\dots,t_i,$ with a birational morphism $\pi\col\mathcal C\ra\mathcal Y,$ which is an isomorphism away from the special fiber. Let $Y\subset\mathcal Y$ be the central fiber. Thus $\pi\col\mathcal C\ra \mathcal Y$ contracts $F_{i+1},\dots,F_\tau$ to $\tau-i$ tacnodes of $Y,$ while $F_1,\dots,F_i$ are  elliptic components of $Y.$ Denote by $Y':=\overline{Y-\cup_{1\le h\le j}F_h}.$

\[
\begin{xy} <16pt,0pt>:
(0,0)*{\scriptstyle}="a"; 
"a"+(-10,3);"a"+(-4,2)**\crv{(-4,0.5)};
"a"+(-10,2);"a"+(-8.5,1.5)**\crv{(-8,4.5)}; 
"a"+(-4,2);"a"+(-3,3)**\crv{(-4,3)}; 
"a"+(-3,3);"a"+(-2,2)**\crv{(-2,3)}; 
"a"+(-2,2);"a"+(-1,3)**\crv{(-2,3)};
"a"+(-2.3,2.2)*{t_3}; 
"a"+(-2,2);"a"+(-2,1.8)**\crv{(-2,2)}; 
"a"+(-2,1.8);"a"+(-2,1)**\crv{(-2.8,1)}; 
"a"+(-2,1.8);"a"+(-2,1)**\crv{(-1.2,1)}; 
"a"+(-6.8,1);"a"+(-5.2,0.7)**\crv{(-6,4)}; 
"a"+(-6,3.3)*{Y};
"a"+(-10.3,2)*{F_1}; 
"a"+(-6.2,-0.3)*{\text{The curve $Y$ for $\tau=3$, $i=2j=2$}}; 
"a"+(-7.1,0.9)*{F_2};
\end{xy}
\begin{xy} <16pt,0pt>:
(0,0)*{\scriptstyle}="a"; 
"a"+(-10,3);"a"+(-4,2)**\crv{(-4,0.5)};
"a"+(-6.8,1);"a"+(-5.2,0.7)**\crv{(-6,4)}; 
"a"+(-7.1,0.9)*{F_2};
"a"+(-4,2);"a"+(-3,3)**\crv{(-4,3)}; 
"a"+(-3,3);"a"+(-2,2)**\crv{(-2,3)}; 
"a"+(-2,2);"a"+(-1,3)**\crv{(-2,3)};
"a"+(-2.3,2.2)*{t_3}; 
"a"+(-2,2);"a"+(-2,1.8)**\crv{(-2,2)}; 
"a"+(-2,1.8);"a"+(-2,1)**\crv{(-2.8,1)}; 
"a"+(-2,1.8);"a"+(-2,1)**\crv{(-1.2,1)}; 
"a"+(-6,3.3)*{Y'};
"a"+(-6.5,-0.3)*{\text{The curve $Y'$ for $\tau=3$, $i=2j=2$}}; 
\end{xy}
\]

\smallskip

We claim the existence of a set $\mathcal T_\xi$ of $2^{\tau-i+j}$ line bundles of $Y$ with $h^0(Y,P)=1$, for every $P\in\mathcal T_\xi,$ such that $\mathcal S_\xi=\{(\pi|_C)^* P \;|\; P\in\mathcal T_\xi\}.$ 

Indeed, recall again that it follows from the proof of \cite[Theorem (2.22)]{Ha} and \cite[Section 2e (a)]{Ha}, that if $Z$ is a curve, $Z'$ is the normalization at a tacnode of $Z$ with the points $p,q$ over the tacnode and $M$ a theta characteristic of $Z$, then there are 2 theta characteristics $M_1,M_2$ of $Z$ whose pull-back to $Z'$ is $M(p+q)$ and such that $h^0(M)+1\equiv h^0(M_i)\;(2)$. 

Now, consider the theta characteristic $G|_{W^\nu\cup F_{j+1}\dots\cup F_i}$ of $W^\nu\cup F_{j+1}\dots\cup F_i$. Since the stable spin curve $\xi=(X,G,\alpha)$ is odd and the restriction of $G$ is an even theta characteristics on $F_1,\dots,F_j$ and an odd theta characteristic on $F_{i+1},\dots,F_\tau,$ by Section \ref{sec2.1} we have that $G|_{W^\nu\cup F_{j+1}\dots\cup F_i}$ is odd (even) if and only if $\tau-i$ is even (odd). Thus we get a set $\mathcal T'_\xi$ of $2^{\tau-i}$ odd theta characteristics of $Y'$ whose pull-back to $W^\nu\cup F_{j+1}\dots\cup F_i$ is $G|_{W^\nu\cup F_{j+1}\dots\cup F_i}(\sum_{i<h\le\tau}(n_{h1}+n_{h2}))$. Now, $W$ is general, thus Lemma \ref{secvan} implies that $h^0(Y', P')=1$ for $P'\in\mathcal T'_\xi.$ Consider the Cartier divisor $D_{\mathcal Y}:=\sum_{1\le h\le j}F_h$ of the total space $\mathcal Y$ of $g\col\mathcal Y\ra B'.$ For $P'\in\mathcal T'_\xi,$ construct the set $\mathcal T_{P'}$ of $2^j$ line bundles of $Y$ by gluing:
\begin{equation}\label{Y'}
\textstyle
P'(\sum_{1\le h\le j}(n_{h1}+n_{h2}))\in\Pic(Y')\,\,\,\,\,\text{and}\,\,\,\,\,G|_{F_h}\in\Pic(F_h)\,\,\,1\le h\le j
\end{equation}

so that, for any $P\in\mathcal T_{P'},$ we have $P^{\otimes 2}\simeq\omega_g(D_{\mathcal Y})\otimes\mathcal O_Y.$ Since $G|_{F_h}$ is non-effective for $1\le h\le j$, it follows that $h^0(Y,P)=1$ for every $P\in\mathcal T_{P'}.$ We show that the set $\mathcal T_\xi:=\cup_{P'\in\mathcal T'_\xi}\mathcal T_{P'}$ of $2^{\tau-i+j}$ line bundles of $Y$ is as required by the claim. Indeed for every $P\in\mathcal T_\xi,$ pick the line bundle $\mathcal P\in\Pic(\mathcal Y)$ extending $P$ to a family of theta characteristics away from the special fiber.
As for Lemma \ref{redcal} (ii), we get:
$$(\pi^*\mathcal P)^{\otimes 2}=\pi^*(\omega_g(D_{\mathcal Y}))=\omega_f(D).$$

Thus $(\pi|_C)^*P\in\mathcal S_\xi.$ We get a map $\mathcal T_\xi\ra\mathcal S_\xi,$ sending $P$ to $(\pi|_C)^*P.$ If $(\pi|_C)^*P_1=(\pi|_C)^*P_2,$ then it follows from the unicity of the extensions in Section \ref{sec2.2} that:
$$(\pi|_C)^*P_1=(\pi^*\mathcal P_1)|_C=(\pi^*\mathcal P_2)|_C=(\pi|_C)^*P_2\Rightarrow\pi^*\mathcal P_1\simeq\pi^*\mathcal P_2.$$

Now, $\pi$ is an isomorphism away from the special fiber and the degree of the restrictions of $P_1$ and $P_2$ to the components of $Y$ are equal, then:
$$(\mathcal P_1)^*\simeq(\pi^*\mathcal P_1)^*\simeq(\pi^*\mathcal P_2)^*\simeq(\mathcal P_2)^*\Rightarrow P_1\simeq P_2.$$

Thus $\mathcal T_\xi\ra\mathcal S_\xi$ is an injection, hence a bijection because the two sets have the same cardinality. It follows that $\mathcal S_\xi=\{(\pi|_C)^* P \;|\; P\in\mathcal T_\xi\}$ and the claim follows.

\smallskip
 
Pick $L=\pi^*P\in\mathcal S_\xi$ for $P\in\mathcal T_\xi$. If $\mathcal P$ is the extension of $P$, then $\pi^*\mathcal P\otimes \mathcal O_C=L.$ Since $L$ has a unique $\pi^*\mathcal P$-smoothable section, also $P$ has a unique $\mathcal P$-smoothable section, hence the unique section $s_P$ of $P$ is $\mathcal P$-smoothable. It follows that $\pi^*s_P$ is limit of fiberwise sections of $\pi^*\mathcal P,$ i.e. $\pi^*s_P$ is the $\pi^*\mathcal P$-smoothable section of $L.$ 
If $P\in\mathcal T_{P'},$ it follows from (\ref{Y'}) that $s_P\in H^0(Y,P)$  
vanishes on $F_1,\dots, F_j$ and restricts to the section of $H^0(Y',P')$ away from $F_1,\dots, F_j.$

Therefore, from Lemma \ref{secvan} we have that:  

\begin{itemize}
\item[(i)]
$\pi^*s_P$ identically vanishes on $F_h$ for $1\le h\le j$;
\item[(ii)]
$\pi^*s_P$ vanishes on a smooth point of $Y$ contained in $F_h$ for $j<h\le i$ ($P$ has degree $1$ on each $F_h$);
\item[(iii)]
$\pi^*s_P$ does not vanish on each curve $F_h$ for $i<h\le\tau$;
\item[(iv)]
$\pi^*s_P$ vanishes on a set $\{l_1,\dots,l_{g-i-j-1}\}_{P'}$ of smooth points of $C$ lying on $W^\nu$ and depending only on $P'$ (the degree of $P$ on the partial normalization of $W$ contained in $Y$ is $g-i+j-1$ and the singular points of $Y$ where $P$ vanishes are the $2j$ points $F_h\cap F_h^c$ for $1\le h\le j$).
\end{itemize}
 
Let $T_{t_h}W$ be the tacnodal tangent to $W$ at $t_h.$ It follows from (i), (ii), (iii), (iv), that, for $P\in \mathcal T_\xi,$ we have: 
\begin{equation}\label{spin-image}
\mu_D(C,\pi^*P)=\spa\{T_{t_1}W,\dots,T_{t_j}W, t_{j+1},\dots,t_i, \spa\{l_1,\dots,l_{g-i-j-1}\}_P\}.
\end{equation}

Notice that $\mu_D(C,\pi^*P)$ is a theta hyperplane of type $(i, j).$ 
 
\smallskip

\noindent
\emph{Third Step: the computation of the multiplicities}. 
By construction, the $2^j$ sections of $\{s_P \;|\; P\in\mathcal T_{P'}\}$ are all equal. Conversely, if $P_1\in\mathcal T_{P'_1}\subset\mathcal T_\xi$ and $P_2\in\mathcal T_{P'_2}\subset\mathcal T_\xi$ with $P'_1\ne P'_2,$ then the set of zeroes of $s_{P_1}$ and $s_{P_2}$ are different away from 
$F_1\cup\dots\cup F_j$, because otherwise $P'_1=P'_2$ ($s_{P_1},s_{P_2}$ do not vanish on components of $Y'$). Thus the set $\{\pi^*s_P \;|\; P\in\mathcal T_\xi\}$ has $2^j 2^{\tau-i}$ sections, each one of which appears $2^j$ times.

Now we vary $\xi$ among odd spin curves supported in $X.$ Consider all the possible odd stable spin curves $(X,G,\alpha)$, where $G|_{F_h}$ varies among the even theta characteristics of $F_h$ for $1\le h\le j$, $G|_{F_h}=\mathcal O_{F_h}$ for $i<h\le\tau$ and $G|_{W^\nu\cup F_{j+1}\cdots\cup F_i}$ is a fixed theta characteristic of $W^\nu\cup F_{j+1}\cdots\cup F_i$.

Each $F_h$ has $3$ even theta characteristics, hence we get a set $\{\xi_1,\dots\xi_{3^j}\}$ of $3^j$ odd stable spin curve. By (i), (ii), (iii), (iv), we have $\pi^*s_{P_1}=\pi^*s_{P_2}$ if and only if $P_1|_{Y'}=P_2|_{Y'}.$ 
 Thus the set $\{\pi^*s_P \;|\; P\in\cup_{1\le r\le 3^j} \mathcal T_{\xi_r}\}$ has $6^j2^{\tau-i}$ sections, each one of which appears $6^j$ times.

Pick the group $\mathcal A_h=\{id,\gamma^h_1,\gamma^h_2,\gamma^h_3\}$ of automorphisms of $F_h$ of Lemma \ref{aut-exc}, for $j<h\le i.$ Since these automorphisms fix each point of $F_h\cap F^c_h,$ they extend to automorphisms both of $C$ and $Y.$ For every $\xi_r=(X,G_r,\alpha_r)\in\{\xi_1,\dots,\xi_{3^j}\},$ define: $$\mathcal U^r:=\{(X,\gamma^*G_r,\gamma^*(\alpha_r)) \;|\; \gamma=\gamma_{j+1}\circ\dots\circ\gamma_i, \gamma_h\in\mathcal A_h\}.$$ 

Notice that by varying $r,$ we get all distinct spin curves. Furthermore, if we set $\mathcal H:=\cup_{1\le r\le 3^j}\{\pi^*s_P \;|\; P\in\cup_{\xi\in\mathcal U^r}\mathcal T_\xi\},$ then we have:
\begin{equation*}
\textstyle
\mathcal H=\{\pi^*\gamma^*s_P \;|\; P\in\cup_{1< r\le 3^j}\mathcal T_{\xi_r}\,;\,\gamma=\gamma_{j+1}\circ\dots\circ\gamma_i, \gamma_h\in\mathcal A_h\}.
\end{equation*}

Indeed if $P\in\mathcal T_{\xi_r},$ then $\gamma^*s_P$ is the unique section of $\gamma^*P,$ and hence it is smoothable. From Lemma \ref{aut-exc}, two sections $\pi^*s_{P_1},\pi^*s_{P_2}$ of $\mathcal H$ are equal away from $\cup_{j<h\le i}{F_h}$ if and only if $\pi^*s_{P_1}=\pi^*\gamma^*s_{P_2}$ where $\gamma$ is an automorphism of $Y$ given by $\gamma=\gamma_{j+1}\circ\dots\circ\gamma_i,$   for $\gamma_h\in\mathcal A_h.$
Thus:

\smallskip

\begin{itemize}
\item[(v)]
the set $\mathcal H$ has $4^{i-j} 6^j 2^{\tau-i}$ sections each one of which vanishes in one point of $F_h,$ for $j<h\le i,$ and is equal, away from $\cup_{j<h\le i}F_h,$ to exactly $4^{i-j} 6^j$ other sections of $\mathcal H.$ 
\end{itemize}

\smallskip

Pick $\pi^*s_P\in\mathcal H.$ By (\ref{spin-image}) we know that $\mu_D(C, \pi^*P)$ is a theta hyperplane of type $(i, j)$ and it follows from (v) that $\mu_D$ sends exactly $4^{i-j} 6^j$ twisted spin curves to $\mu_D(C,\pi^*P).$ 
By Example \ref{banana}, each twisted spin curve has multiplicity 1 in the central fiber of $S^\nu_f\ra B',$ thus the following claim implies that a theta hyperplane of type $(i, j)$ has multiplicity $4^{i-j} 6^j.$

\smallskip

Claim: $\psi(C, \pi^*\overline{P})\ne\mu_D(C, \pi^*P)$ if $\pi^*\overline{P}\notin\mathcal H$.

Indeed, pick $\xi_r=(X,G_r,\alpha_r)\in\{\xi_1\dots\xi_{3^j}\}.$ Let $\overline{\xi}$ be another odd spin curve of $C$ supported on $X.$ If $\mu_D(C,\pi^*P)=\psi(C,\pi^*\overline{P})$ for $P\in\mathcal S_{\xi_r}$ and $\overline{P}\in\mathcal S_{\overline\xi},$ then $\pi^*s_P$ and $\pi^*s_{\overline{P}}$ vanish on one point of $F_h$ for $j< h\le i,$ and they are equal away from $\cup_{j<h\le i}F_h.$ Thus $\overline{\xi}\in\mathcal U^r$ for some $r.$

Let $\xi'=(X', G', \alpha'),$ where $X'\ne X.$ If $P'\in\mathcal T_{\xi'}$ and $P\in\mathcal T_\xi,$ then $\psi(C,\pi^*P')\ne\mu_D(C,\pi^*P)$ because the type of the two hyperplanes is different. The tacnodal case of the Theorem is done.

\smallskip

\noindent
Assume now that the singular points of $W$ are exactly $\gamma$ cusps $c_1,\dots,c_\gamma.$ Let $\mathcal W\ra B$ be a general smoothing of $W$ to theta-generic curves and $f\col\mathcal C\ra B'$ be its stable reduction as in Lemma \ref{redcal}, with central fiber $C.$ Set $n_h:=F_h\cap F^c_h.$ Pick the divisor of $\mathcal C$ given by $D=\sum_{1\le h\le\gamma}F_h$ and the smooth curve $S^-_{\mathcal N_D},$ all of whose points over $0\in B'$ are supported on $C.$ Consider:
\[
\SelectTips{cm}{11}
\begin{xy} <16pt,0pt>:
\xymatrix{
 S^\nu_f \UseTips\ar[r]^{\nu_D}\ar[dr]_{\psi}& S^-_{\mathcal N_D} \ar[d]^{\mu_D}\\
 & J_{\mathcal W'}
}
\end{xy}
\]

such that $\mu_D\circ\nu_D=\psi.$ Denote by $\nu_0\col S_f^\nu\ra S^-_{\omega_f}$ the normalization of $S^-_{\omega_f}.$

Fix a stable odd spin curve $\xi=(X,G,\alpha)$ of $C,$ a point of $S^-_{\omega_f}.$ It is supported on the blow-up $X$ of $C$ at all of its nodes. Assume that:
\begin{equation}
G|_{F_h}\text{ is even for }1\le h\le k;\,\,G|_{F_h}=\mathcal O_{F_h} \text{ for }k<h\le\gamma.
\end{equation}

Of course, $G|_{W^\nu}$ is a theta characteristic of $W^\nu$.

Consider the $D$-twisted spin curve $(C, L)$, where $L$ satisfies $L^{\otimes 2}\simeq \omega_C\otimes\mathcal O_f(D)$ and has restrictions: 
\begin{equation*}
\textstyle
L|_{F_h}=G|_{F_h} \text{ for }1\le h\le k;\,\,\,\,L|_{F_h}=G|_{F_h}=\mathcal O_{F_h} \text{ for }k<h\le\gamma.
\end{equation*}
\begin{equation*}
\textstyle
L|_{W^\nu}=G|_{W^\nu}(\sum_{1\le h\le\gamma}n_h).
\end{equation*}

Arguing as in the tacnodal case, we see that there exists a representative $(X, G, \alpha)$ of $\xi$ such that $G$ and $L$ are limits of the same family of theta characteristics. Thus $\nu_0^{-1}(\xi)=(C,L)\in S^\nu_f$ because, being $C$ of compact type, then $S^-_{\omega_f}\ra B'$ is \'etale and the cardinality of $\nu_0^{-1}(\xi)$ is 1. In order to describe the morphism $\psi\col S^\nu_f\ra J_{\mathcal W'},$ it suffices to find the images of the $D$-twisted spin curves via the morphism $\mu_D\col S^-_{\mathcal N_D}\ra J_{\mathcal W'}.$
Arguing as in the tacnodal case, one can show that, if $s$ is the smoothable section of $L$, then $s$ identically vanishes on $F_1\cup\dots\cup F_k,$ $s$ does not vanishes on $F_{k+1}\cup\dots\cup F_\gamma,$ $s$ has a set $\{l_1,\dots,l_{g-k-1}\}_L$ of $g-k-1$ zeroes  on smooth points of $C$ on $W^\nu$ and two different sections have different sets of zeroes. We have: $$\mu_D(C,L)=\spa\{c_1,\dots,c_k,\spa\{l_1,\dots l_{g-k-1}\}_L\},$$ 

hence $\mu_D(C,L)$ is of type $k.$ If we change the $3^k$ even theta characteristics of $F_1,\dots,F_k,$ then $\mu_D(C,L)$ does not change. Now, $(C, L)$ has multiplicity 1 in the central fiber of $S^\nu_f,$ then a theta hyperplane of type $k$ has multiplicity $3^k.$

\smallskip

\noindent
The case of a curve with tacnodes and cusps follows by repeating word by word the proofs of the case of a curve with just tacnodes and of a curve with just cusps.
\end{proof}

The technique used to prove Theorem \ref{mult} applies also to nodal curves with at most two components, as is shown in \cite[4.1.1.]{Pa}.

\begin{Exa}
The reader can check that $N_g=\sum 4^{i-j} 6^j 3^k t^j_{ik}(W)$. For example, if $W$ has 1 tacnode and 1 cusp, we have $\widetilde{g}=g-3$ and from Theorem \ref{hyper}:  $$t^0_{00}=2 t^1_{11}=2 N_{g-3};\,\,\,\,t^0_{01}=2 t^1_{10}=2 N^+_{g-3};\,\,\,\,t^0_{10}=t^0_{11}=2^{2(g-3)}.$$ As expected: $$\sum 4^{i-j} 6^j 3^k t^j_{ik}(W)=36N_{g-3}+28N^+_{g-3}=32\cdot 2^{2(g-3)}-4(N^+_{g-3}-N_{g-3})=N_g.$$
\end{Exa}

\subsection{Spin curves over non-stable curves}\label{sec5.4} We can conclude with a geometric meaningful definition of spin curves over non-stable curves.

\begin{Def}
Let $W$ be an irreducible curve, whose singularities are cusps and tacnodes. 
A \emph{spin curve of $W$}
is a triple $(C,T,L),$ where:
\begin{itemize} 
\item[(i)]
$C$ is the central fiber of the stable reduction $f\col\mathcal C\ra B'$ of a general smoothing of $W$;
\item[(ii)]
$T=\mathcal O_f(D)\otimes\mathcal O_C$ is a twister of $C,$ where $D$ is a  divisor of $\mathcal C$ given by the sum with coefficient $1$ of all the elliptic tails lying over the cusps of $W$ and of some elliptic tails lying over the tacnodes of $W;$
\item[(iii)]
$L\in\Pic C$ is a square root of $\omega_C\otimes T.$
\end{itemize}
\end{Def}

The curve $S^\nu_f$ of (\ref{nonst-spin}) has a description in terms of spin curves of $W,$ thus the limits of odd theta characteristics degenerating to $W$ are spin curves of $W$.

\subsection*{Acknowledgments}
I thank L. Caporaso for her warm support and her valuable guidance, 
C. Casagrande, G. Mondello, E. Esteves, E. Sernesi, D. Testa for fundamental suggestions and the anonymous referee for his detailed report and his very constructive comments.

\end{document}